\documentclass[12pt,halfline,a4paper]{ouparticle}
\graphicspath{{Figure/}}
\usepackage{multirow}
\usepackage{multicol}
\usepackage{amssymb}
\usepackage{amsmath}
\usepackage{bbm}
\usepackage{ntheorem}
\usepackage{amsbsy}
\usepackage[round]{natbib}   
\usepackage{authblk}
\newtheorem{theorem}{Theorem}
\newtheorem{assumption}{Assumption} 

\newtheorem{lemma}{Lemma}
\newtheorem{remark}{Remark}
 \def\bx{{\mathbf x}}
 \def \bb{{\pmb{\beta}}}
 \def\b1{{\mathbf 1}}
\def\bI{{\mathbf I}} 
\def\bz{{\pmb{z}}} \def \bnu{{ \pmb{\nu}}}
\def \bv{\pmb{\nu}}
\def\br{{\mathbf r}}
\def\RR{\mathbb{R}} 

\def \mL{ \mathcal{L}}
\def \bo{{\mathcal{O}}}
\def \b0{{\bf 0}}
\def\bD{{\mathbf D}}
\def \blp{\left (}
\def \brp{\right )}

\newcommand{\argmin}{\operatornamewithlimits{argmin}}
\newcommand{\indep}{\perp\mkern-9.5mu\perp}
\newcommand*\diff{\mathop{}\!\mathrm{d}}
\newcommand*{\QEDB}{\hfill\ensuremath{\square}}
\newtheorem{example}{Example}[section]
\begin{document}

\title{High-dimensional robust approximated $M$-estimators for mean regression with asymmetric data}
\author[1]{Bin Luo \thanks{Corresponding author. Email: bin.luo2@duke.edu}} 
\author[2]{Xiaoli Gao \thanks{Gao is partially supported by Simons Foundation Grant: SF359337 and SF854940.}}
\affil[1]{Department of Biostatistics and Bioinformatics, Duke University}
\affil[2]{Department of Mathematics and Statistics, The University of North Carolina at Greensboro}

\abstract{
Asymmetry along with heteroscedasticity or contamination often occurs with the growth of data dimensionality. In ultra-high dimensional data analysis, such irregular settings are usually overlooked for both theoretical and computational convenience. In this paper, we establish a framework for estimation in high-dimensional regression models using Penalized Robust Approximated quadratic M-estimators (PRAM). This framework allows general settings such as random errors lack symmetry and homogeneity, or covariates are not sub-Gaussian. To reduce the possible bias caused by data's irregularity in mean regression, PRAM adopts a loss function with an adaptive robustification parameter. Theoretically, we first show that, in the ultra-high dimension setting, PRAM estimators have local estimation consistency at the minimax rate enjoyed by the LS-Lasso. Then we show that PRAM with an appropriate non-convex penalty in fact agrees with the local oracle solution, and thus obtain its oracle property. Computationally, we compare the performances of six PRAM estimators (Huber, Tukey's biweight or Cauchy loss function combined with Lasso or MCP penalty function). Our simulation studies and real data analysis demonstrate satisfactory finite sample performances of the PRAM estimator under different irregular settings.}

 \date{}


\maketitle

\section{Introduction}
\label{sec1}
Asymmetry along with heteroscedasticity or contamination often occurs with the growth of data dimensionality. In high-dimensional settings, particularly when random errors follow irregular distributions such as asymmetry and heteroscedasticity, simultaneous mean estimation and variable selection are still of interest in many applications. For example, in economics where asymmetric data is prevalent, it is of interest to study how mean GDP is affected by many features. The natural asymmetries and high skewness of the GDP data are also addressed in \cite{hess1997asymmetric} and \cite{dew2021skewness}. One suggestion is to use the data transformation such as the log transformation for improving the detection of important variables \citep{lutkepohl2012role}. However,  we may be still interested in building a mean regression model on the original data. 

In this paper, we are interested in high-dimensional mean regression that is robust to the following irregular settings: (a) the data are not symmetric due to the skewness of random errors (\cite{fan2017estimation}); (b) the data are heteroscedastic (\cite{daye2012high}, \cite{wang2012quantile}); and (c) the data are contaminated in both response and a large number of variables (\cite{rousseeuw2005robust}). However, the above irregular settings are often overlooked for high-dimensional data analysis, especially for theoretical development.
 
Penalized M-estimator have been studied in high-dimensional robust regression (e.g., \cite{huber1964robust}, \cite{lambert2011robust}, \cite{gao2010asymptotic}, \cite{wang2013l1},  \cite{loh2017statistical}, \cite{loh2021scale}). However,  most of the previous work does not estimate the conditional mean function, and/or require the error distribution to be symmetric and/or homogeneous. Recently, \cite{fan2017estimation} proposed a so-called RA-Lasso estimator to estimate the mean function under asymmetric random errors. In their work, the Huber loss with a diverging parameter is adopted to reduce the bias when the error distribution is asymmetric. \cite{fan2017estimation} obtained nice asymptotic properties of the RA-Lasso estimator, and proved its estimation consistency at the minimax rate enjoyed by LS-Lasso. The RA-Lasso estimator was also investigated in \cite{sun2019adaptive} when random errors are required to have bounded $k$th moment with $k>1$. \cite{pan2021iteratively} further extended the RA-Lasso to nonconvex regularization by applying the iteratively reweighted L1-penalized algorithm and established the oracle properties.  

However, the Huber loss approximation used in \cite{fan2017estimation}, \cite{sun2019adaptive} and \cite{pan2021iteratively} does not downweight the very large residual due to its non-decreasing $\Psi$-function. \cite{shevlyakov2008redescending} showed that M-estimators given by non-decreasing $\Psi$-function do not possess finite variance sensitivity. It means that the asymptotic variance can be largely affected if the assumed model is only approximately true. In that paper, the authors proposed to consider re-descending M-estimators with $\Psi$-function re-descending to zero to address this problem. They further showed that a re-descending M-estimator can be designed by maximizing the minimum variance sensitivity under a global minimax criterion. For instance, the Smith's estimator and Tukey's biweight estimator are two optimal M-estimators with minimax variance sensitivity for a class of densities with a bounded variance and a bounded fourth moment, respectively \citep{shevlyakov2008redescending}. Therefore it is tempting to consider a re-descending $M$-estimator in complex high-dimensional settings with heavy data contamination.

For decades both the theoretical and computational results in penalized re-descending M-estimator in high-dimensional settings have been very limited, due to the non-convex loss functions. Recently, \cite{loh2017statistical} established a form of local statistical consistency for the high-dimensional $M$-estimators allowing both the loss and  penalty functions to be non-convex. However, this study does not address the problem of asymmetry and heteroscedasticity, and therefore is not applicable for mean regression. Also, their numerical studies neglect settings for asymmetric data and lack comparisons among different $M$-estimations.

In this paper, we consider high-dimensional mean regression in more general irregular settings: the data can be contaminated or include possible large outliers in both random errors and covariates, the random errors may lack symmetry and homogeneity. In particular, we investigate both statistical and computational properties of high-dimensional mean regression 
in the penalized $M$-estimator framework with adaptive robustification parameters. This framework allows both the loss function and the penalty to be non-convex. Our perspective is different from \cite{loh2017statistical} since all loss functions considered in our study approximate to the quadratic loss with an adaptive robustification parameter. To be more specific,  we proposed a class of Penalized Robust Approximated quadratic $M$-estimators (PRAM) to obtain a robust mean estimator under irregular settings (a-c) mentioned above. The PRAM uses a family of loss functions with an adaptive parameter to control the robustness as well as the discrepancy to the quadratic loss. With the adaptively chosen robustification parameter, PRAM estimators are able to reduce the bias induced by asymmetric error distribution and meanwhile preserve the robustness to approximate the mean estimators. This approximation was also used in \cite{fan2017estimation}, \cite{sun2019adaptive} and \cite{pan2021iteratively}, where the Huber loss was used to obtain some robustness.
     
Our theoretical contributions in this paper include the investigation of statistical properties for a class of PRAM estimators with only weak assumptions on both random errors and covariates. In particular, we first provide upper bounds of the estimation errors and introduce  sufficient conditions under which  a PRAM estimator has local estimation consistency at the same rate as the minimax rate enjoyed by the LS-Lasso. We then  show that the PRAM estimator actually equals the local oracle solution with the correct support if an appropriate non-convex penalty is used. Based on this oracle result we further establish the asymptotic normality of the PRAM estimators. By choosing the robustification parameter adaptively to sample size, dimension and moments, our theoretical results are applicable for a wide class of PRAM estimators that are robust to general irregular settings, when the dimensionality of data grows with the sample size at an almost exponential rate.

Computationally, we implement the PRAM estimator through a two-step optimization procedure: obtain an initial PRAM estimator using the convex Huber loss function with the Lasso penalty in the first step, then optimize the program of the desired PRAM estimator using this initial in the second step. We also investigate the performance of six PRAM estimators generated from three types of loss function approximation (the Huber loss, Tukey's biweight loss and Cauchy loss) combined with two types of penalty functions (the Lasso and MCP penalties). While our numerical results demonstrate satisfactory finite sample performance of the PRAM estimators under general irregular settings, it suggests that in practice, when the data are heavy-tailed or contaminated, a well-behaved PRAM estimator can be chosen by considering a re-descending loss function approximation and a concave penalty, using the RA-Lasso as an initial.   


\paragraph{Related Works:} we end this section by highlighting a few things on how our work is different from some recent related work on mean estimation:

As introduced earlier, the Huber loss approximation considered in   \cite{fan2017estimation}, \cite{sun2019adaptive} and \cite{pan2021iteratively} waives the symmetry requirement by allowing the parameter of Huber loss to diverge. The idea is that by controlling the divergence rate of the parameter, while preserving certain robustness, the Huber loss becomes `closer' to the $\ell_2$ loss and thus potentially reduces the bias when the error distribution is asymmetric. Our work in this paper relaxes the convexity restriction of loss functions and answers the question on how in general a loss function with strong robustness should approximate to the $\ell_2$ loss to achieve the estimation consistency at the minimax rate. In addition, our framework also allows concave penalties and therefore inherits certain oracle properties under some conditions. 

 While the work in \cite{fan2017estimation} and  \cite{pan2021iteratively} required the existence of at least second moment of random errors conditional on covariates, \cite{sun2019adaptive} relaxed this condition and studied the RA-Lasso requiring the existence of the ($1+\delta$)-th moment  for $\delta>0$.  They investigated the phase transition between $\delta \ge 1$ and $0<\delta<1$ and showed the estimation consistency of RA-Lasso at a slower rate for the latter. Additionally, \cite{sun2019adaptive} relaxed the sub-Gaussian or sub-exponential assumption on covariates by applying a truncation approach to the covariates.  Their assumption requires the existence of bounded uniform fourth moment of the covariate vector. The truncation approach is also applied to both the response and covariates in \cite{fan2021shrinkage}.  
In this paper, we also only require the existence of the $k$-th moment ($k>1$) of random errors for PRAM. Besides the estimation consistency, we further establish the oracle property for only $k>1$ in general. We additionally relax the assumption of covariates by incorporating weight functions in the extension of PRAM estimators. We later show that we even place no moment requirement on the covariates with a proper choice of the weight functions.

    

 The remainder of our paper is organized as follows.  In Section 2, we introduce the basic setup regarding PRAM estimators and corresponding generalizations. In Section 3, we establish the local estimation consistency for the PRAM estimators under sufficient conditions.  For non-convex regularized PRAM estimators, we also present our statistical theory concerning the selection consistency and the asymptotic normality of PRAM estimators.  We discuss the implementation of PRAM estimators including both the computational algorithm and the tuning parameter selection in Section 4. In section 5, we conduct some simulation studies to demonstrate the performance of the PRAM estimators under different settings. We also apply those PRAM estimators for NCI-60 data analysis and  illustrate all results in Section 6. Section 7 concludes and summarizes the paper. All technical proofs are relegated to the Appendix.

\paragraph{Notation:} We use bold symbols to denote matrices or vectors. For a matrix or a vector $\bv$, we write $\bv^T$ to denote its transpose. We write $\|\cdot\|_1$ and  $\|\cdot\|_2$ to denote the $L_1$ norm and the $L_2$ norm of a vector, respectively. For a function $g: \RR^p \mapsto \RR$, we write $\nabla g$ to denote a gradient of the function. We write $u_+$ to denote $\max(u,0)$ for any $u \in \RR$. We write $f \precsim g$ for functions $f$ and $g$ when $f \le cg$ for some universal constant $c >0$. We write $f \asymp g$ to mean that $f \precsim g$ and $g \precsim f$ hold simultaneously.

\section{The PRAM method}\label{sec2}
\subsection{Model settings}\label{sec2-1}
Consider an ultra high-dimensional linear regression model
\begin{equation} \label{linear model}
   	 y_i=\bx_i^T\bb^*+\epsilon_i,
\end{equation}
where $\bx_i=(x_{i1},\cdots,x_{ip})^T$ for $i=1,\cdots,n$ are independent and identically distributed (i.i.d) $p$-dimensional covariate vectors such that $E(\bx_i)=\bf 0$, $\{\epsilon_i\}_{i=1}^n$
 are independent errors such that $E(\epsilon_i \mid \bx_i)=0$ and thus we allow the conditional heteroscedasticity. Note $\bb^*=(\beta^*_1,\cdots,\beta^*_p)^T\in\RR^p$  with $p \gg n$ is an $s$-sparse conditional mean coefficient vector (only includes $s$  nonzero elements).
 

Our model settings permit the existence of all following irregular settings on both $\epsilon_i$s and $\bx_i$s: (a) asymmetry of $\epsilon_i$; (b) heteroscedasticity of $\epsilon_i$ and possible dependence of $\epsilon_i$ on $\bx_i$; (c) data contamination of $\epsilon_i$ and $\bx_i$.

We are interested in penalized mean regression estimators such that 
 \begin{equation} \label{eq:estimator}
 \hat{\bb}\in \argmin_{\|\bb\|_1\le R} \left\{ \mL_{\alpha,n}(\bb)+\rho_\lambda(\bb)\right\},
 \end{equation}
where $\mL_{\alpha, n}$ is the empirical loss function and $\rho_\lambda$ is a penalty function which encourages the sparsity in the solution. Here $\alpha >0$ is an adaptive robustification parameter controlling the robustness, which is allowed to diverge. We consider the loss function $\mL_{\alpha, n}$ satisfying
\begin{equation} \label{population condition}
    \lim_{\alpha \to \infty} E[\nabla \mL_{\alpha, n}(\bb^*)] =\b0.
\end{equation}
 This condition in (\ref{population condition}) relaxes the condition, $E[\nabla \mL_{\alpha, n}(\bb^*)] =\b0$ for each $\alpha >0$, required in \cite{loh2017statistical}, which is not guaranteed with the lack of homogeneity and symmetry. Condition (\ref{population condition}) permits the random error to be  heterogeneous and/or asymmetric, as long as $E[\nabla \mL_{\alpha, n}(\bb^*)]$ converges to $\b0$ with a diverging $\alpha$.

We also include the side condition $\|\bb\|_1 \le R$ for $R \ge \|\bb^*\|_1$ in the penalized optimization problem in (\ref{eq:estimator}) to guarantee the existence of local/global optima, when the penalized objective function is non-convex. In real applications, we can choose $R$ to be a sufficiently large number.


\subsection{Penalty functions}\label{sec2-2}
Since the coefficients vector $\bb^*$ is assumed to be $s$-sparse in  the high-dimensional linear regression model in \eqref{linear model}, we only consider penalties that generate sparse solutions. 
In particular, we require the penalty function $\rho_{\lambda}$ in \eqref{eq:estimator} to satisfy following properties listed in Assumption 1. 
\begin{assumption}[Penalty Function Assumptions] \label{as:penalty}
	The penalty function is coordinate-separable such that $\rho_\lambda(\bb)=\sum_{j=1}^{p} \rho_\lambda(\beta_j)$ for some scalar function $\rho_\lambda: \RR \mapsto \RR$. In addition,
\begin{enumerate}
\item[(i)] the function $t \mapsto \rho_\lambda(t)$ is symmetric around zero and $\rho_\lambda(0)=0$;
\item[(ii)] the function $t \mapsto \rho_\lambda(t)$  is non-decreasing on $\RR^+$;
\item[(iii)] the function $t \mapsto \frac{\rho_\lambda(t)}{t}$ is non-increasing on $\RR^+$;
\item[(iv)] the function $t \mapsto \rho_\lambda(t)$  is differentiable for $t \ne 0$;
\item[(v)] $ \lim_{t \to 0^+} \rho_\lambda'(t) = \lambda$;
\item[(vi)] there exists $\mu > 0$ such that the function $t \mapsto \rho_{\lambda} (t) + \frac{\mu}{2}t^2$ is convex;
\item[(vii)] there exists $\delta \in (0, \infty)$ such that $\rho_{\lambda}'(t)=0$ for all $t \ge \delta \lambda$. 
\end{enumerate} 
\end{assumption}
Those properties in Assumption 1 are related to the penalty functions  studied in  \cite{loh2013regularized} and \cite{loh2017statistical}, where $\rho_{\lambda}$ is said to be $\mu$-amenable if $\rho_{\lambda}$ satisfies conditions (i)-(vi) for $\mu$ defined in (vi). If $\rho_{\lambda}$ also satisfies condition (vii), we say that $\rho_{\lambda}$ is $(\mu, \delta)$-amenable. Some popular choices of amenable penalty functions include Lasso \citep{tibshirani1996regression}, SCAD \citep{fan2001variable}, and MCP  \citep{zhang2010nearly}. We only present the expression of Lasso and MCP as follows since they will be used in the numerical studies.
\begin{itemize}
    \item The Lasso penalty, $\rho_{\lambda}(t)=\lambda |t|$, is $0$-amenable but not $(0,\delta)$-amenable for any $\delta <\infty$.\\
\item The  MCP penalty,
\begin{equation*}
\rho_{\lambda}(t)={\rm sign}(t) \lambda \int_{0}^{|t|}\left(1-\frac{z}{\lambda b}\right)_+ dz,
\end{equation*}
where $b>0$ is a fixed parameter. The MCP penalty is also $(\mu, \delta)$-amenable with $\mu=\frac{1}{b}$ and $\delta=b$.   
\end{itemize}
It has been shown that the folded concave penalty, such as SCAD or MCP, possesses better variable selection properties than the convex penalty like the Lasso.
\subsection{Loss functions} \label{sec2-3}
From the linear model setting in Section \ref{sec2-1}, we know $E(y_i | \bx_i)=\bx_i^T\bb^*$. We are interested in finding a well-behaved mean regression estimator of $\bb^*$. Since we consider a general setting discussed in Section \ref{sec2-1}, we wish to study the empirical loss function $\mL_{\alpha,n}$ that are robust to outliers and/or heavy-tailed distribution. Let $\ell_\alpha: \RR \mapsto \RR $ denote a residual function, or a loss function, defined on each observation pair ($\bx_i, y_i$). The corresponding empirical loss function for (\ref{eq:estimator}) is then given by
\begin{equation} \label{eq:M-estimator}
\mL_{\alpha,n}(\bb) = \frac{1}{n} \sum_{i=1}^{n}\ell_\alpha(y_i - \bx_i^T\bb).
\end{equation}
With a well-chosen non-quadratic function $\ell_\alpha$, the penalized mean regression estimator from (\ref{eq:estimator}) can be robust to outliers or heavy-tailed distribution in the additive noise term $\epsilon_i$. However, such a mean estimator can be biased when the conditional distribution of $\epsilon_i$ is not symmetric. 

To reduce the bias, we consider a family of loss functions with an adaptive robustification parameter satisfying (\ref{population condition}) to  approximate the traditional quadratic loss. In particular, we require the following approximation:
\begin{equation} \label{eq:approximate loss}
\lim_{\alpha \to \infty} \ell_\alpha(u)=\frac{1}{2} u^2, \quad \forall u \in \RR.
\end{equation}
The empirical loss function satisfying (\ref{eq:approximate loss}) is called a robust approximated quadratic loss function. The following approximations take robust loss functions such as the Huber loss, Tukey's biweight loss and Cauchy loss to  approximate the quadratic loss function: 
\begin{itemize}
    \item 
 Huber Approximation
\begin{equation*}
\ell_\alpha(u)=\begin{cases}
\frac{u^2}{2} & \text{if } |u| \le \alpha, \\
\alpha |u| - \frac{\alpha^2}{2} & \text{if } |u| \ge \alpha;
\end{cases} 
\end{equation*}
\item Tukey's biweight Approximation
\begin{equation*}
\ell_\alpha(u)=\begin{cases}
\frac{\alpha^2}{6}(1-(1-\frac{u^2}{\alpha^2})^3) & \text{if } |u| \le \alpha, \\
\frac{\alpha^2}{6} & \text{if } |u| \ge \alpha;
\end{cases}
\end{equation*}
\item  Cauchy Approximation
\begin{equation*}
\ell_\alpha(u)=\frac{\alpha^2}{2}\log(1+\frac{u^2}{\alpha^2}).
\end{equation*}
\end{itemize}
It is straightforward to verify that all above three loss functions satisfy equation (\ref{eq:approximate loss}). In addition, the Tukey's biweight loss and Cauchy loss produce re-descending $M$-estimators. In the robust regression literature, we call an $M$-estimator re-descending if there exists $u_0>0$ such that $|\ell_\alpha'(u)|=0$ or decrease to 0 smoothly, for all $|u| \ge u_0$. In that case, large residuals can be downweighted. See more discussions in \cite{muller2004redescending, shevlyakov2008redescending}.

\subsection{PRAM estimators and the extensions} \label{sec2-4}
A class of PRAM estimators takes the form:
 \begin{equation} \label{eq:PRAM-estimator}
 \hat{\bb} \in \argmin_{\|\bb\|_1\le R} \left\{\frac{1}{n} \sum_{i=1}^{n}\ell_\alpha(y_i-\bx_i^T \bb)+\rho_\lambda(\bb)\right\},
 \end{equation}
where the penalty function $\rho_\lambda$ satisfies Assumption \ref{as:penalty} and the loss function $\ell_\alpha$ is a scalar function satisfying the Approximation Equation \eqref{eq:approximate loss}.

Whereas a PRAM estimator in equation (\ref{eq:PRAM-estimator}) takes into account the contamination or heavy-tailed distribution in asymmetric additive error, a single outlier in $\bx_i$ may still cause the corresponding estimator to perform arbitrarily badly. By downweighting large values of $\bx_i$, we extend the class of PRAM estimators to  
 \begin{equation} \label{eq:GPRAM-estimator}
 \hat{\bb} \in \argmin_{\|\bb\|_1\le R} \left\{
 \frac{1}{n} \sum_{i=1}^{n} \frac{w(\bx_i)}{v(\bx_i)}\ell_\alpha((y_i-\bx_i^T \bb)v(\bx_i))+\rho_\lambda(\bb)\right
 \},
 \end{equation}
where $w,v$ are weight functions mapping from $\RR^p$ to $\RR^+$. When $w \equiv v\equiv 1$, (\ref{eq:GPRAM-estimator}) is reduced to the PRAM defined in (\ref{eq:PRAM-estimator}).  Some examples of the weight functions $w$ and $v$ are as follows.
\begin{itemize}
    \item Mallows weight \cite{mallows1975some}
    \begin{equation} \label{eq:mallow}
        w(\bx)=\min\left\{1,\frac{b} {\|B\bx\|_2}\right\} \ \text{or}\ \min\left\{1,\frac{b^2} {\|B\bx\|_2^2}\right\}, v(\bx)\equiv 1,
    \end{equation}
    for some $b>0$ and $B \in \RR^{p \times p}$. Note that $\|w(\bx)\bx\|_2$ is bounded for fixed values of $b$ and $B$. 
    
    \item
    Hill-Ryan weight\cite{hill1977robust}
        \begin{equation} \label{eq:hill}
        v(\bx)\equiv w(\bx),
    \end{equation}
    where $w(\bx)$ is defined such that $\|w(\bx)\bx\|_2$ is bounded (e.g., (\ref{eq:mallow})).
    
    \item Schweppe weight \cite{merrill1971bad}
    \begin{equation*} 
            w(\bx)=\frac{1}{\|Bx\|_2}, v(\bx)=\frac{1}{w(\bx)}.
    \end{equation*}
    
\end{itemize}
The above $w(\bx)$s are able to bound the influence function of the M-estimator \cite{loh2017statistical, hampel1986robust} and downweight the contribution of high leverage points in (\ref{eq:GPRAM-estimator}). In addition, the $v(x)$ in Hill-Ryan's method further shrinks the residuals based on the leverage weight of $\bx$. For the Schweppe weight, if $\ell_\alpha'$ is re-descending to zero, the $v(\bx)$ additionally dampens the effect of a leverage point by inflating the residual according to the leverage of $\bx$.

For the rest of the paper, we specify the PRAM estimator with the Huber approximation, Tukey's biweight approximation and Cauchy approximation as the HA-type, TA-type and CA-type PRAM estimator, respectively. In particular, we also specify a PRAM estimator using a re-descending loss function approximation (e.g., Tukey's biweight approximation and Cauchy approximation) a re-descending PRAM estimator. Additionally, we classify a PRAM estimator with the Lasso penalty and MCP penalty as the Lasso-type and MCP-type PRAM estimator correspondingly.

\section{Statistical Properties}
\label{sec3}


\subsection{Estimation Consistency} \label{sec3-1}
As in (\ref{eq:GPRAM-estimator}), we consider a class of PRAM estimators with the loss function,
\begin{equation} \label{eq:gloss}
\mL_{\alpha,n}(\bb)=\frac{1}{n} \sum_{i=1}^{n} \frac{w(\bx_i)}{v(\bx_i)}\ell_\alpha((y_i-\bx_i^T \bb)v(\bx_i)).
\end{equation}
To obtain the estimation consistency, we make following assumptions on  $\ell_\alpha$.
\begin{assumption}[Loss Function Assumptions] \label{as:loss}
$\ell_\alpha: \RR \mapsto \RR$ is a  scalar function for $\alpha > 0$ with the existence of the first derivative $\ell_\alpha'$ everywhere and the second derivative $\ell_\alpha''$ almost everywhere. In addition,
	\begin{enumerate} 
		\item[(i)] there exists a constant $0 < k_1 < \infty$ such that $|\ell_\alpha'(u)| \le k_1 \alpha$ for all $ u \in \RR$;
		\item[(ii)] for all $\alpha >0$, $\ell_\alpha'(0)=0$ and $\ell_\alpha'$ is Lipschitz such that $|\ell_\alpha'(x) - \ell_\alpha'(y)| \le k_2|x-y|$ for all $x, y \in \RR$ and some $0 <k_2 <\infty$;
		\item[(iii)] for some $k > 1$, there exists a constant $d_1>0$ such that $|1-\ell_\alpha''(u)| \le d_1|u|^k\alpha^{-k}$ for almost all $|u| \le \alpha$.
	\end{enumerate}
\end{assumption}


Note that Assumption \ref{as:loss} (i) indicates that the magnitude of $\ell_\alpha'$ is bounded from above at the same rate of $\alpha$ so that the PRAM estimator can achieve robustness. Assumption \ref{as:loss} (ii) implies $|\ell_\alpha'(u)| \le k_2|u|$ for all $u \in \RR$ and $|\ell_\alpha''(u)| \le k_2 $ for almost every $u \in \RR$. In particular, the loss functions we study in this paper actually satisfy Assumption \ref{as:loss} (ii) with  $k_2=1$, showing that $\ell_\alpha$ is bounded by the quadratic loss function $u^2/2$ for any $\alpha$. Assumption \ref{as:loss} (iii) indicates that for almost all $u \in \RR$, $\ell''_\alpha$ converges point-wisely to 1 with at least the order of $\alpha^{-k}$ for $k\ge 1$. 

The above assumptions cover a wide range of loss functions, including the Huber loss, Hampel loss, Tukey's biweight loss and Cauchy loss.

\begin{remark}
By some simple math, we can show that $\lim_{\alpha \to \infty} \ell'_{\alpha}(u)= u$ for all $u\in \RR$ based on Assumption \ref{as:loss}. Suppose in addition that $\ell_\alpha(0)=0$, we can further obtain the approximation  in (\ref{eq:approximate loss}), indicating that Assumption 2 alone gives sufficient conditions for $\ell_\alpha$ to approximate the quadratic loss.
\end{remark}
\begin{remark}
By the dominated convergence theorem, we have
\begin{equation*}
\begin{split}
    \lim_{\alpha \to \infty} E[\nabla \mL_{\alpha, n}(\bb^*)] &= \lim_{\alpha \to \infty} E[w(\bx_i)\bx_i\ell_\alpha'(\epsilon_i v(\bx_i))] = E[w(\bx_i)\bx_i (\epsilon_i v(\bx_i))]\\
    &= E[w(\bx_i)\bx_iE(\epsilon_i\mid \bx_i) v(\bx_i))]=\b0 .
\end{split}
\end{equation*}
Thus Condition (\ref{population condition}) in Section \ref{sec2-1} is guaranteed, and therefore, heterogeneous and/or asymmetric error is allowed.
\end{remark}

We now make some weak assumptions on both random error $\epsilon$ and covariate vector $\bx$ for the investigation of the approximation error. 
\begin{assumption}[Error and Covariate Assumptions] \label{as:x} For $w(\bx)$ and $v(\bx)$ given in \eqref{eq:GPRAM-estimator}, 
the random error $\epsilon$ with $E[\epsilon \mid \bx]=0$ and covariate vector $\bx$ with $E[\bx]=\b0$ satisfy:
	\begin{enumerate}
		\item[(i)] $E[E(|\epsilon|^k \mid \bx)v(\bx)^k]^2 \le M_{k} < \infty$, for $k > 1$ in Assumption 2(iii);
		\item[(ii)] $\sup_{\|u\|_2=1} E[v(\bx)\bx^Tu]^{2k}=q_k < \infty$, for $k >1 $ in Assumption 2(iii);
		\item[(iii)] $0<k_l<\lambda_{min}(E[w(\bx)v(\bx)\bx\bx^T])$ and $\lambda_{max}(E[w(\bx)^2\bx\bx^T]) < k_u$;
		\item[(iv)] for any $\pmb \nu \in \RR^p$, $w(\bx)\bx^T\pmb \nu$ is sub-Gaussian with parameter at most $k_0^2\|\pmb \nu\|_2^2$.
	\end{enumerate}
\end{assumption}

Note that condition (i) requires only the existence of $k$-th conditional moment of $\epsilon$ with $k>1$, indicating that this condition is independent of the distribution of $\epsilon$ itself and can hold for many heavy-tailed or skewed distributions. If $w(\bx) \equiv v(\bx)\equiv 1$, conditions (ii) and (iv) hold when $\bx_i^T\pmb \nu$ is sub-Gaussian for any $\pmb \nu \in \RR^p$. In this case, Assumption \ref{as:x} becomes conditions (C1-C3) in \cite{fan2017estimation}, except that we only require $k>1$ instead of $k \ge 2$. If covariate $\bx$ is contaminated or heavy-tailed distributed, conditions (ii)-(iv) nonetheless hold with some proper choices of $w(\bx)$ and $v(\bx)$ (e.g., (\ref{eq:mallow}) and (\ref{eq:hill})). Thus, the sub-Gaussian assumption on $\bx$ is potentially relaxed.

\begin{remark}
In the relevant work of \cite{sun2019adaptive} and \cite{fan2021shrinkage}, a prepossessing truncation approach is applied to $\bx$ for heavy-tailed data. Their assumption requires that $\bx$ has bounded uniform fourth moment. In contrast,  our weighting methods downweight the contribution of the leverage points in (\ref{eq:gloss}) and are able to further relax the assumption on $\bx$. Specifically, if the Mallows weight in (\ref{eq:mallow}) is used, condition (ii) only require the bounded uniform ($2+\delta$)th moment of $\bx$ for any $\delta>0$. If the Hill-Ryan weight in (\ref{eq:hill}) is adopted, we place no moment requirement on $\bx$ as $\|v(\bx)\bx\|_2$ is bounded.

\end{remark}

Let $\bb_{\alpha}^*$ be a local non-penalized population minimizer under the PRAM loss,
\begin{equation}\label{eq:beta_alpha_star}
\bb_{\alpha}^* \in \argmin_{\|\bb - \bb^*\|_2 \le R_0} \left \{ E \left[\frac{w(\bx)}{v(\bx)}\ell_\alpha((y-\bx^T\bb)v(\bx))\right] \right \},
\end{equation}
for some $0 < R_0< \infty$. Note that $\bb_{\alpha}^*$ is a local minimizer of (\ref{eq:beta_alpha_star}) within a neighborhood of $\bb^*$. The restriction $\|\bb - \bb^*\|_2 \le R_0$ is included to guarantee the existence of local optima, when the penalized objective function is non-convex.
If the regularization parameter $\lambda$ in equation (\ref{eq:GPRAM-estimator}) converges to 0 sufficiently fast, then $\hat{\bb}$ is a natural unpenalized $M$-estimator of $\bb_{\alpha}^*$ for any $\alpha>0$. Whereas $\bb_{\alpha}^*$ differs from $\bb^*$ in general, $\bb_{\alpha}^*$ is expected to converge to $\bb^*$ when $\alpha\to \infty$, due to the approximation in (\ref{eq:approximate loss}) for PRAM. The rate of the approximation error $\|\bb_{\alpha}^* - \bb^*\|_2$ is established in the following Theorem 1. 

\begin{theorem} \label{Tm:1}
Under the Assumption \ref{as:loss} and \ref{as:x}, there exists a universal positive constant $C_1$, such that $\|\bb_{\alpha}^*-\bb^*\|_2 \le d{\alpha}^{1-k}$, where $d=2^k C_1k_l^{-1}\sqrt{k_u}(\sqrt{M_k}+R_0^k\sqrt{q_k}){\alpha}^{1-k}$. Here $k$, $k_l$, $k_u$, $M_k$, $q_k$ appear in Assumption \ref{as:loss} and \ref{as:x}, and $R_0$ appears in (\ref{eq:beta_alpha_star}).
\end{theorem}
Theorem \ref{Tm:1} gives an upper bound of the approximation error between the true parameter vector and the non-penalized PRAM population minimizer. The approximation error vanishes when the adaptive robustification parameter $\alpha\to \infty$. It vanishes faster if a higher moment of $\epsilon | \bx$ exists. In fact, Theorem 1 demonstrates that the approximation of $\ell_\alpha$ to the quadratic loss helps to reduce the bias induced by the asymmetry on $\epsilon$. If we let $\ell_\alpha$ in (\ref{eq:gloss}) be the Huber loss and $w(\bx) \equiv v(\bx)\equiv 1$, Theorem 1 gives the upper bound of the approximation error studied in \cite{fan2017estimation}.

In order to obtain the estimation  consistency for the PRAM estimator in (\ref{eq:GPRAM-estimator}), we also require the loss function $\mL_{\alpha,n}$ to satisfy the following uniform Restricted Strong Convexity (RSC) condition. 

\begin{assumption}[Uniform RSC condition]\label{as:RSC}
	There exist $\gamma$, $\tau$, $\alpha_0>0$ and a radius $r>0$ such that for all $\alpha \ge \alpha_0$, the loss function
	$\mL_{\alpha,n}$ in  \eqref{eq:GPRAM-estimator} satisfies
	\begin{equation} \label{eq:URSC}
	\langle \nabla \mL_{\alpha,n}(\bb_1)- \nabla \mL_{\alpha,n}(\bb_2), \bb_1-\bb_2 \rangle \ge \gamma \|\bb_1-\bb_2\|_2^2 - \tau \frac{\log p}{n}\|\bb_1 - \bb_2\|_1^2,
	\end{equation}
	where $\bb_j \in \RR^p$ such that $\|\bb_j-\bb^*\|_2\le r$ for $j \in \{1, 2\}$.
\end{assumption}

Note that the uniform RSC assumption is only imposed
on $\mL_{\alpha,n}$ inside the ball of radius $r$ centered at $\bb^*$. Thus the loss function used for robust regression can be wildly nonconvex while it is away from  the origin. The radius $r$ essentially specifies a local ball centered around $\bb^*$ in which stationary points of the PRAM estimator are well-behaved.

\begin{remark}
In  \cite{loh2013regularized} and \cite{loh2017statistical}, the RSC condition were imposed on a specific loss function. Although  Assumption \ref{as:RSC} requires that the RSC condition is satisfied uniformly over a family of loss functions generated from a range of $\alpha$, this assumption is in fact not stronger: Assumption \ref{as:RSC} holds naturally if there exists $\alpha_0>0$ such that $\mL_{\alpha_0,n}$ satisfies Assumption 2 and inequality (\ref{eq:URSC}) for some $\gamma, \tau >0$. We further establish the uniform RSC condition in the Appendix .
\end{remark}

We present our main estimation consistency results on the PRAM estimator in the following Theorem \ref{Tm: est}.

\begin{theorem}\label{Tm: est} 
Suppose $\mL_{\alpha,n}$   in \eqref{eq:GPRAM-estimator} satisfies Assumption \ref{as:loss} and the uniform RSC condition in Assumption \ref{as:RSC}. Suppose $\rho_\lambda$ is $\mu$-amenable with $\frac{3}{4}\mu < \gamma$ in Assumption \ref{as:penalty}.
Assume the random error and covariates satisfy Assumption \ref{as:x}. Let $\max\{(\frac{2d}{R_0})^{\frac{1}{k-1}},C_{21}(\frac{n}{\log p})^{\frac{\zeta}{k-1}}\} \le \alpha \le C_{31} (\frac{n}{\log p})^{1-\zeta}$ and $\lambda \ge C_5(\frac{\log p}{n})^{\zeta}$ with $0 < \zeta \le \min\{\frac{k-1}{k}, \frac{1}{2}\}$. Let $\hat{\bb}$ be a local PRAM estimator in  the uniform RSC region. Then for $n > C_0r^{-\frac{1}{\zeta}}s^{\frac{1}{2\zeta}}\log p$, with probability at least $1-2\exp\{-C_4 \log p\}$, $\hat{\bb}$ exists and satisfies
	\begin{equation*}
	   \|\hat{\bb}-\bb^*\|_2 \le \frac{24C_5}{4\gamma-3\mu}\sqrt{s} \left(\frac{\log p}{n}\right)^\zeta  \text{ and  } \|\hat{\bb}-\bb^*\|_1 \le \frac{96C_5}{4\gamma-3\mu}s\left (\frac{\log p}{n} \right)^\zeta.
	\end{equation*}
\end{theorem}
The nonasymptotic result of Theorem \ref{Tm: est} holds even when the random errors lack symmetry and homogeneity, and the regressors lack sub-Gaussian assumption. It provides upper bounds of estimation errors with high probability, if the adaptive robustification parameter  $\alpha$ diverges at the rate $(\frac{n}{\log p})^{\frac{\zeta}{k-1}} \precsim \alpha \precsim (\frac{n}{\log p})^{1-\zeta}$. In particular, if $\epsilon|\bx$ has a finite variance ($k \ge 2$),  Theorem \ref{Tm: est} implies $ \|\hat{\bb}-\bb^*\|_2 \precsim \sqrt{\frac{s\log p}{n}}$ and $\|\hat{\bb}-\bb^*\|_1 \precsim \sqrt{\frac{s^2\log p}{n}}$ with $\lambda \asymp \sqrt{\frac{\log p}{n}}$. Hence, Theorem \ref{Tm: est} guarantees that a PRAM estimator within the local region of radius $r$ performs as well as the LS-Lasso under the sub-Gaussian random error.  The rate range of $\alpha$ stated in Theorem \ref{Tm: est} requires that $\alpha$ should diverge faster enough, e.g., faster than $\bo \left ((\frac{n}{\log p})^{\frac{\zeta}{k-1}} \right )$, to reduce the bias sufficiently but meanwhile not too fast, e.g., slower than $\bo \left ((\frac{n}{\log p})^{1-\zeta} \right)$, in order to preserve certain robustness of a PRAM estimator. The existence of a higher moment of $\epsilon|\bx$ (a larger $k$) actually allows $\alpha$ to diverge at a lower rate.  

\begin{remark}
The convergence rate of the PRAM estimator relies on the divergence rate of the adaptive robustification $\alpha$. If we let $\alpha$ diverge with $\zeta=\min(\frac{k-1}{k}, \frac{1}{2})$, we obtain the optimal rate of convergence as described in Theorem 3 in \cite{sun2019adaptive}  for the adaptive Huber estimator with the Lasso penalty. While \cite{sun2019adaptive} requires a fixed divergence rate $\alpha \asymp (\frac{n}{\log p})^{1-\zeta}$, we allows $\alpha$ to  diverge at a flexible range of rate as shown in Theorem \ref{Tm: est}. The range indicates a trade off between $\alpha$, $k$ and $\zeta$: a lower convergence rate (a smaller $\zeta$) or a higher moment of $\epsilon|\bx$ (a larger k) leads to a broader effective range of $\alpha$ for $\hat{\bb}$ to have the above non-asymptotic result. As we will show in the next subsection, $\alpha$ is only allowed to diverge with $\zeta <\min(\frac{k-1}{k}, \frac{1}{2})$ for establishing the oracle properties, if $\epsilon|\bx$ has  at most a finite variance, i.e., $1<k\le2$, and $s$ grows with $n$.
\end{remark}

\begin{remark}
The proof of Theorem \ref{Tm: est} in  Appendix B reveals that the estimation consistency result also holds for all local stationary points of the penalized objective function  in  (\ref{eq:estimator}) in Section \ref{sec2-1}. Here $\tilde{\bb}$ is a stationary point of the optimization in (\ref{eq:estimator}) if 

\begin{equation*}
    \langle \nabla \mL_{\alpha, n}(\tilde{\bb}) + \nabla \rho_{\lambda}(\tilde{\bb}), \bb - \tilde{\bb} \rangle \ge 0,
\end{equation*}
for all feasible $\bb$ in a neighbour of $\tilde{\bb}$.
Note that stationary points include both the interior local maxima as well as all local and global minima.
Hence, Theorem \ref{Tm: est} guarantees that all stationary points within the ball of radius $r$ centered at $\bb^*$ have local statistical consistency at the minimax rate enjoyed by the LS-Lasso.
\end{remark}

\subsection{Oracle Properties}
In this section, we establish the oracle properties for the PRAM estimators  in program (\ref{eq:GPRAM-estimator}). We first define a local oracle estimator $\hat{\bb}^\mathcal{O}=(\hat{\bb}_S^{\mathcal{O}}, \b0_{S^c})$ with
\begin{equation} \label{eq:oracle}
    \hat{\bb}_{S}^{\mathcal{O}}= \argmin_{\bb \in \RR^S:\|\bb - \bb^*\|_2 \le r} \left\{ \mL_{\alpha, n} (\bb) 
        \right \},
\end{equation}
where $S=\{j: \beta^*_j \ne 0 \}$ and the side condition $\|\bb - \bb^*\|_2 \le r$ is used to guarantee that the oracle estimator $\hat{\bb}^\mathcal{O}$ lies within the uniform RSC region. 
Let  $\beta^*_{\min}= \min_{j \in S} |\beta^*_j|$ denotes a minimum signal strength on $\bb^*$. Our oracle result shows that when the penalty $\rho_\lambda$ is $(\mu, \delta)$-amenable, those stationary points of the PRAM estimator in program (\ref{eq:GPRAM-estimator}) within the local neighborhood of $\bb^*$ can be unique and agree with the oracle estimator $\hat{\bb}^{\mathcal{O}}$ with high probability, as stated in the following theorem.

\begin{theorem} \label{Tm: oracle}
Suppose the penalty $\rho_{\lambda}$ is $(\mu, \delta)$-amenable and conditions in Theorem \ref{Tm: est} hold. Suppose in addition that $v(\bx)\bx_j$ is sub-Gaussian for all $j \in \{1, \cdots,p\}$, $\|\bb^*\|_1 \le \frac{R}{2}$ for some $R>\frac{192\lambda s}{4\gamma-3\mu}$. Let $\lambda=C_{71}(\frac{\log p}{n})^{\zeta}$ for $0 < \zeta \le \min\{\frac{k-1}{k}, \frac{1}{2}\}$, $\beta^*_{\min} \ge C_{72} (\frac{\log p}{n})^{\zeta}$, and $n \ge C_{01} s \log p$ for a sufficiently large constant $C_{01}$. Suppose $\alpha$ satisfies $C_{22}s^{\frac{1}{k-1}}(\frac{n}{\log p})^{\frac{\zeta}{k-1}} \le \alpha \le C_{31} (\frac{n}{\log p})^{1-\zeta}$ and $s = \bo \left ((\frac{n}{\log p})^{k(1-\zeta)-1} \right)$. Let $\tilde{\bb}$ be a stationary point of program (\ref{eq:GPRAM-estimator}) in the uniform RSC region. Then we have

\begin{equation*}
    P(\tilde{\bb}=\hat{\bb}^{\mathcal{O}}) \ge 1-C_8 \exp(-C_{41}s^{-\frac{1}{\zeta}}\log p).
\end{equation*}
\end{theorem}
 
Note that the lower bound of $\alpha$ in Theorem \ref{Tm: oracle} is higher than the bound in Theorem \ref{Tm: est}, with a ratio $\bo \left (s^{\frac{1}{k-1}} \right)$. To guarantee the existence of an optimal $\alpha$, $s$ cannot grow with $n$ too fast, provided  $s = \bo \left ((\frac{n}{\log p})^{k(1-\zeta)-1} \right)$ for $k > 1$.
However,  if the number of non-zero parameters $s$ is finite, the condition on $s$ is trivial and the effective range of $\alpha$ becomes the same as one in Theorem \ref{Tm: est}. The feasibility condition $\|\bb^*\| \le \frac{R}{2}$ instead of $R$ in Theorem 2, is for the technical proof.
\begin{remark}
The condition  $s = \bo \left ((\frac{n}{\log p})^{k(1-\zeta)-1} \right)$ implies that the existence of a higher moment of $\epsilon|\bx$ (a larger $k$) allows $s$ to grow faster. Specifically, a diverging $s$ requires $\alpha$ to diverge with $\zeta <\frac{k-1}{k}$. If $k>2$, we can let $\zeta=\frac{1}{2}$. Then the condition of the minimum signal strength becomes $\beta^*_{\min} \ge C_{72} \sqrt{\frac{\log p}{n}}$. If $1 < k \le 2$, it requires $\beta^*_{\min} \ge C_{72} (\frac{\log p}{n})^{\zeta}> C_{72} (\frac{\log p}{n})^{\frac{k-1}{k}} \ge C_{72} \sqrt {\frac{\log p}{n}}$. Hence, a smaller $k$ may requires a larger minimum signal strength to guarantee the oracle result with diverging $s$.

\end{remark}

Applying our Theorem \ref{Tm: oracle} and the standard results for M-estimators with a diverging number of parameters in \cite{he2000parameters}, we can obtain the following theorem concerning the asymptotic normality of any stationary point of the program (\ref{eq:GPRAM-estimator}). For the sake of simplicity, we only provide the result under $w(\bx) \equiv v(\bx)\equiv 1$. The result of a weighted PRAM can be derived accordingly. 

\begin{theorem} \label{Tm: normality}
Suppose conditions in Theorem \ref{Tm: oracle} hold with $k > 2 $ and the loss function $\mL_{\alpha, n}$ given in (\ref{eq:gloss}) is twice differentiable within the $\ell_2$-ball of radius r around $\bb^*$. Suppose for all $\alpha >0$, $\ell_\alpha''$ is Lipschitz such that $|\ell_\alpha''(x) - \ell_\alpha''(y)| \le k_3|x-y|$ for all $x, y \in \RR$ and some $0< k_3 <\infty$. Suppose in addition that $\alpha > (2C_9/k_l)^{1/k}$ and $\alpha ^{1-k} = o(n^{-1/2})$. 
Let $\tilde{\bb}$ be a stationary point of program (\ref{eq:GPRAM-estimator}) in the uniform RSC region. If $\frac{s \log^3 s}{n}  \rightarrow 0$, then $\|\tilde{\bb} - \bb^*\|_2 = \bo_p \left (\sqrt{\frac{s}{n}} \right)$. If $\frac{s^2 \log s}{n} \rightarrow 0$, then for any $\bv \in \RR^p$, we have
\begin{equation*}
    (\sqrt{n}/\sigma_{\bv}) \cdot \bv^T(\tilde{\bb} - \bb^*) \xrightarrow{d} N(0,1),
\end{equation*}
where 
\begin{equation*}
   \sigma_{\bv}^2 = \bv^T_S E[(\nabla^2\mL_{\alpha, n}(\bb^*))_{SS}]^{-1} Var(\ell_\alpha'(\epsilon_i)(\bx_i)_S)  E[(\nabla^2\mL_{\alpha, n}(\bb^*))_{SS}]^{-1} \bv_S.
\end{equation*}
\end{theorem}

The condition $\alpha ^{1-k} = o(n^{-1/2})$ indicates that $\alpha$ should diverge at least faster than $n^{\frac{1}{2(k-1)}}$, in addition to the rate stated in Theorem \ref{Tm: oracle}. Here the condition $\alpha > (2C_9/k_l)^{1/k}$ is required to guarantee the invertibility of matrix $E[(\nabla^2\mL_{\alpha, n}(\bb^*))_{SS}]$.
\begin{remark}
To further understand the condition $\alpha ^{1-k} = o(n^{-1/2})$, we take $\alpha = \bo \left (\sqrt{\frac{n}{\log p}} \right)$ as an example, the fastest divergence rate indicated in Theorem \ref{Tm: oracle}  with $\zeta = \frac{1}{2}$. Then the condition requires $ \log p = o(n^{\frac{k-2}{k-1}})$. Therefore, the asymptotic normality result holds only when a higher moment ($k>2$) of $\epsilon|\bx$ exists. 
\end{remark}


\section{Implementation of the PRAM estimators} \label{sec:4}

Note that the optimization in (\ref{eq:estimator}) may not be a convex optimization problem since we allow both the loss function $\mL_{\alpha,n}$ and $\rho_\lambda$ to be non-convex. To obtain the corresponding stationary point,  we use the composite gradient descend algorithm \citep{nesterov2013gradient}. Denoting $q_\lambda(\bb) = \lambda\|\bb\|_1-\rho_\lambda(\bb)$ and $\bar{L}_{\alpha,n}(\bb)=\mL_{\alpha,n}(\bb)-q_\lambda(\bb)$, we can rewrite the program as
\begin{equation*}
\hat{\bb} \in \argmin_{\|\bb\|_1 \le R} \left\{ \bar{L}_{\alpha,n}(\bb) + \lambda\|\bb\|\right\}.
\end{equation*}
Then the composition gradient iteration is given by
\begin{equation} \label{eq:interate}
\bb^{t+1} \in \argmin_{\|\bb\|_1 \le R} \left\{ \frac{1}{2} \| \bb - (\bb^t -  \eta \nabla \bar{L}_{\alpha,n}(\bb^t))\|_2^2 + \eta\lambda \|\bb\|_1  \right\},
\end{equation}
where $\eta >0$ is the step size for the update and can be determined by the backtracking line search method described in \cite{nesterov2013gradient}. A simple calculation shows that the iteration in (\ref{eq:interate}) takes the form
\begin{equation*}
\bb^{t+1} = S_{\eta\lambda } \left(\bb^t -  \eta\nabla \bar{L}_{\alpha,n}(\bb^t)\right),
\end{equation*}
where $S_{\eta \lambda}(\cdot)$ is the soft-thresholding operator defined as 
\begin{equation*}
[ S_{\eta\lambda }(\bb)]_j =  {\rm sign}(\beta_j) \left(|\beta_j| - \eta\lambda\right)_+.
\end{equation*}
We further adopt the two-step procedure discussed in \cite{loh2017statistical} to guarantee the convergence to a stationary point for the non-convex optimization problem:
\begin{enumerate}
    \item[1.] Run the composite gradient descent using the convex Huber loss function with the convex Lasso penalty to get an initial PRAM estimator.
     \item[2.] Run the composite gradient descent on the desired high-dimensional PRAM estimator using the initial PRAM estimator from step 1.
\end{enumerate}

For tuning parameters selection, the optimal values of $\alpha$ and $\lambda$ are chosen by a two-dimensional grid search using cross-validation. In particular, the searching grid is formed by partitioning a rectangle uniformly in the scale of ($\alpha,  \log(\lambda)$). In order to obtain the robust selection of tuning parameters, their optimal values are found by the combination that minimizes the cross-validated trimmed mean squared prediction error, where $t_{cv}\%$ of squared prediction errors are trimmed from each end before the mean is computed. In our numerical study, we choose $t_{cv}=10$. 

\section{Simulation Studies} \label{sec:simulation}
In this section, we assess the performance of the PRAM estimators by considering different types of loss and penalty functions through various models. Similar simulation settings are also investigated in \cite{fan2017estimation}. The data is generated from the following model
\begin{equation*} 
   	 y_i=\bx_i^T\bb^*+\epsilon_i.
\end{equation*}
We choose the true regression coefficient vector as $\bb^*=(\textbf{3}_5^T, \textbf{2}_5^T, \textbf{1.5}_5^T, \textbf{0}_{p-15}^T)^T$, where the first 15 elements consist of 5 numbers of $3$, $2$, $1.5$ receptively and the rest are 0. In all simulation settings, we let $\textit{n=100}$ and $\textit{p=500}$.

\begin{example} \label{ex:1} {\bf (Homogeneous case)}
The covariates vector $\bx_i$s are generated from a multivariate normal distribution with mean $\bf 0$ and covariance  $\bI_p$ independently.  The random errors $\epsilon_i=e_i-E[e_i]$, where $e_i$ are generated independently from the following 5 scenarios:
\begin{enumerate}
\item[(a).] $N(0,4)$: Normal with mean 0 and variance 4;
\item[(b).] $\sqrt{2}t_3$: $\sqrt{2}$ times the $t$-distribution with degrees of freedom $3$;
\item[(c).] MixN: Equal mixture of Normal distributions N(-1, 4) and N(8, 1);
\item[(d).]  LogNormal: Log-normal distribution such that $ e_i = \exp(1.3z_i)$, where $z_i\sim N(0,1)$.
\item[(e).] Weibull: Weibull distribution with the shape parameter $ 0.3$ and the scale parameter $0.15$.
\end{enumerate}
 \end{example}


We consider three types of loss functions equipped with adaptive robustification parameters (the Huber, Tukey's biweight and Cauchy losses) and two types of penalty functions (the Lasso and MCP penalties). Thus it produces 6 different PRAM estimators: HA-Lasso, TA-Lasso, CA-Lasso, HA-MCP, TA-MCP and CA-MCP.  Note the HA-Lasso becomes the RA-Lasso estimator in \cite{fan2017estimation}, which has been demonstrated to perform better than the Lasso and R-Lasso, especially when the errors were asymmetric and heavy-tailed (LogNormal and Weibull). Thus in our simulation we skip those comparisons and only evaluate the performance of all above 6 PRAM estimators. Their performances on both mean estimation  and variable selection  under the five scenarios were reported by the following five measurements:
\begin{enumerate}
\item[(1)] $L_2$ error, which is defined as $\|\hat{\bb} - \bb^*\|_2$.
\item[(2)] $L_1$ error, which is defined as $\|\hat{\bb} - \bb^*\|_1$.
\item[(3)] Model size (MS), the average number of selected covariates.
\item[(4)] False positives rate (FPR), the percent of selected but unimportant covariates:
\begin{equation} \label{eq:FPR}
    FPR=\frac{|\hat{S} \bigcap S^c|}{|S^c|} \times 100\%.
\end{equation}
\item[(5)] False negatives rate (FNR), the percent of non-selected but important covariates:
\begin{equation}
    FNR=\frac{|\hat{S}^c \bigcap S|}{|S|} \times 100\%.
\end{equation}
\end{enumerate}
Here $\hat{S}=\{j: \hat{\beta}_j \ne 0 \}$ and $S=\{j: \beta^*_j \ne 0 \}$.
The model considered in Example \ref{ex:1} is homogeneous, in which the error distribution is independent of covariate $\bx$. We also assess the performance of PRAM estimators under heteroscedastic model in the next example.

\begin{example} \label{ex:2} {\bf (Heteroscedastic case)} We generate the data from 
\begin{equation*}
y_i=\bx_i^T\bb^*+c^{-1}(\bx_i^T \bb^*)^2 \epsilon_i,
\end{equation*}
where the constant $c=\sqrt{3}\|\bb^*\|_2^2$ makes $E[c^{-1}(\bx_i^T \bb^*)^2]^2=1$.
We also consider $\bx_i \sim N(\b0, \bI_p)$ and generate the random error $\epsilon$ from the same five scenarios described in Example \ref{ex:1}.
\end{example}

Finally, we design a simulation setting to evaluate the performance of the generalized PRAM estimators under non-Gaussian covariates.
\begin{example} \label{ex:3} {\bf (Non-Gaussian $\bx$ case)} Similar to Example \ref{ex:1}, except that the covariate $\bx$ in $20\%$ of observations are first generated from independent chi-square variables with 10 degrees of freedom, and then recentered to have mean zero.
\end{example}

For all three examples described above, we run 100 simulations for each scenario. In Example \ref{ex:3}, we consider the generalized PRAM estimators with $v(\bx) \equiv 1$ and  $w(\bx)=\min\left\{1,\frac{4} {\|\bx\|_{2}}\right\}$ (Mallows weight in equation (8)) or  $w(\bx)=\min\left\{1,\frac{4} {\|\bx\|_{\infty}}\right\}$. Here both $w(x)$s can downweight the contribution of high leverage points as they tend to have the larger $\|\bx\|_2$ or $\|\bx\|_{\infty}$.   For all six PRAM estimators, tuning parameters $\lambda$ and $\alpha$ are chosen optimally by 10-fold cross-validation, with $\alpha$ ranges in ($0.1\sqrt{\frac{n}{\log p}}, 10\sqrt{\frac{n}{\log p}}$) and $\lambda$ ranges in ($0.01\sqrt{\frac{\log p}{n}}, 2.5\sqrt{\frac{\log p}{n}}$). These ranges are motivated from Theorem 2. The mean values out of $100$ iterations (with the standard deviation in parentheses) are reported in Table \ref{tb:homo}, \ref{tb:heter}, \ref{tb:contaminatex_inf} and \ref{tb:contaminatex_mallow}, respectively.

We have two findings based on results in Table \ref{tb:homo} and \ref{tb:heter}. Firstly, all the MCP-type PRAM estimators largely outperform the Lasso-type estimators in all the measurements, rendering satisfactory finite sample performances under different settings. This is consistent with the oracle property of the PRAM estimators using a proper non-convex penalty stated in Theorem 3. Secondly, for estimators with the same penalty, although all estimators perform comparably for light-tailed settings ($N(0,4)$ and MixN), the TA-type and CA-type PRAM estimators outperform the HA-type estimators using the same penalty in heavy-tailed settings ($\sqrt{2}t_3$, LogNormal and Weibull). This is actually not surprising due to the following two facts: (1) re-descending M-estimators can achieve the minimax variance sensitivity under certain global minimax criterion \citep{shevlyakov2008redescending};  (2) the HA-Lasso estimation is used as the initial in the optimization process of TA-type and CA-type PRAM estimators. Note that the error terms $c^{-1}(\bx_i^T \bb^*)^2 \epsilon_i$ in the heteroscedastic model have the same variance as those in the homogeneous model, however, their distribution possesses heavier tails. Hence in the heteroscedastic model, except for a few errors being far away on the tail, most of the others get even closer to the center. This fact explains why the performances in Table \ref{tb:heter} are consistently better than those in Table \ref{tb:homo}. 


\begin{table}
  \caption{\label{tb:homo}Simulation results under the homogeneous model with standard normal covariates in Example \ref{ex:1}. The mean $L_2$ error, $L_1$ error, MS, FPR (\%) and FNR (\%) out of 100 iterations are displayed. Standard deviations are listed in parentheses.}
 \centering
 \resizebox{\textwidth}{!}{%
 \fbox{%
\begin{tabular}{clrrrrrr}
                           & \multicolumn{1}{l|}{}                        & HA-Lasso & TA-Lasso & CA-Lasso& HA-MCP  & TA-MCP & CA-MCP \\ \hline
\multirow{3}{*}{N(0,4)}    & \multicolumn{1}{l|}{$L_2$  error} &   3.3 (0.9)           &    3.31 (0.94)         &   3.3 (0.9)             &    0.99 (0.47)         &   1.01 (0.53)          &     0.94 (0.27)         \\
                           & \multicolumn{1}{l|}{$L_1$  error}&   17.75 (4.2)           &    17.79 (4.33)        &   17.72 (4.12)            &    3.29 (2.07)        &   3.34 (2.22)         &     3.06 (1.01)              \\
                           & \multicolumn{1}{l|}{MS}&   67.32 (9.11)           &    66.99 (10.15)        &   66.95 (10.13)            &    17.21 (2.47)        &   16.84 (2.46)         &     16.71 (2.36)              \\
                           & \multicolumn{1}{l|}{FPR, FNR}  &   10.85, 2.13         &   10.8, 2.53      &   10.79, 2.4       &    0.46, 0.27    &   0.4, 0.53   &     0.35, 0.07                        \\ \hline

\multirow{3}{*}{$\sqrt{2}t_3$}    & \multicolumn{1}{l|}{$L_2$  error} &   3.59 (0.96)           &    3.62 (1)         &   3.56 (1)             &    1.18 (0.9)         &   1.13 (0.89)          &     1.14 (0.93)         \\
                           & \multicolumn{1}{l|}{$L_1$  error}&   19.09 (5.03)           &    19.26 (5.14)        &   19.04 (5.12)            &    3.95 (3.65)        &   3.78 (3.56)         &     3.76 (3.65)              \\
                           & \multicolumn{1}{l|}{MS}&   63.72 (9.76)           &    64.07 (11.39)        &   65.13 (9.55)            &    16.85 (2.11)        &   16.7 (2.53)         &     16.2 (2.17)              \\
                           & \multicolumn{1}{l|}{FPR, FNR}  &   10.14, 3         &   10.23, 3.53      &   10.43, 3.13       &    0.43, 1.6    &   0.4, 1.53   &     0.31, 1.87                        \\ \hline

\multirow{3}{*}{MixN}       & \multicolumn{1}{l|}{$L_2$  error} &   3.48 (0.78)           &    3.48 (0.79)         &   3.5 (0.8)             &    1.25 (0.71)         &   1.27 (0.73)          &     1.25 (0.69)         \\
                           & \multicolumn{1}{l|}{$L_1$  error}&   18.99 (3.71)           &    18.99 (3.72)        &   19.05 (3.8)            &    4.2 (2.97)        &   4.17 (2.79)         &     4.11 (2.72)              \\
                           & \multicolumn{1}{l|}{MS}&   68.12 (8.85)           &    68.14 (9.06)        &   67.65 (9.4)            &    17.52 (3.57)        &   17.05 (3.8)         &     17.06 (5.43)              \\
                           & \multicolumn{1}{l|}{FPR, FNR}  &   11, 1.6         &   11.01, 1.73      &   10.92, 2       &    0.55, 0.93    &   0.47, 1.47   &     0.46, 1.2                        \\ \hline

\multirow{3}{*}{LogNormal}   & \multicolumn{1}{l|}{$L_2$  error} &   4.66 (1.2)           &    4.56 (1.13)         &   4.5 (1.24)             &    2.13 (2.05)         &   1.74 (1.68)          &     2.12 (1.97)         \\
                           & \multicolumn{1}{l|}{$L_1$  error}&   23.84 (6.2)           &    23.75 (5.63)        &   23.44 (6.15)            &    7.69 (8.52)        &   5.88 (6.61)         &     7.4 (7.88)              \\
                           & \multicolumn{1}{l|}{MS}&   57.16 (11.44)           &    60.68 (14.11)        &   60.64 (12.13)            &    16.7 (3.61)        &   16.03 (2.69)         &     16.29 (7.03)              \\
                           & \multicolumn{1}{l|}{FPR, FNR}  &   8.97, 8.93         &   9.68, 8.53      &   9.68, 8.73       &    0.62, 8.73    &   0.41, 6.53   &     0.58, 10.13                        \\ \hline

\multirow{3}{*}{Weibull}    & \multicolumn{1}{l|}{$L_2$  error} &   3.91 (1.06)           &    3.63 (1.05)         &   3.46 (1.08)             &    1.35 (1.43)         &   0.94 (1.15)          &     1.03 (1.26)         \\
                           & \multicolumn{1}{l|}{$L_1$  error}&   19.62 (5.38)           &    19.15 (5.36)        &   18.17 (5.65)            &    4.64 (5.73)        &   3.18 (4.5)         &     3.42 (4.69)              \\
                           & \multicolumn{1}{l|}{MS}&   55.37 (11.91)           &    64.12 (11.96)        &   63.5 (8.98)            &    16.15 (2.47)        &   15.65 (1.76)         &     15.44 (1.65)              \\
                           & \multicolumn{1}{l|}{FPR, FNR}  &   8.51, 5.87         &   10.26, 4.13      &   10.09, 3.07       &    0.36, 4.07    &   0.2, 2.13   &     0.18, 2.87 

\end{tabular}}}
\end{table}

\begin{table}
  \caption{\label{tb:heter}Simulation results under the heteroscedastic model with standard normal covariates in Example \ref{ex:2}. The mean $L_2$ error, $L_1$ error, MS, FPR (\%) and FNR (\%) out of 100 iterations are displayed. Standard deviations are listed in parentheses. }
 \centering 
 \resizebox{\textwidth}{!}{%
 \fbox{%
\begin{tabular}{clrrrrrr}
                           & \multicolumn{1}{l|}{}                        & HA-Lasso & TA-Lasso & CA-Lasso& HA-MCP  & TA-MCP & CA-MCP \\ \hline
\multirow{3}{*}{N(0,4)}    & \multicolumn{1}{l|}{$L_2$  error} &   2.84 (0.81)           &    2.94 (0.91)         &   2.72 (0.84)             &    0.55 (0.35)         &   0.55 (0.19)          &     0.6 (0.21)         \\
                           & \multicolumn{1}{l|}{$L_1$  error}&   14.74 (4.14)           &    15.45 (4.84)        &   14.13 (4.38)            &    1.78 (1.16)        &   1.73 (0.64)         &     1.91 (0.82)              \\
                           & \multicolumn{1}{l|}{MS}&   61.56 (9.65)           &    63.25 (11.03)        &   62.11 (8.42)            &    15.68 (1.27)        &   15.28 (0.87)         &     15.43 (1.71)              \\
                           & \multicolumn{1}{l|}{FPR, FNR}  &   9.62, 0.67         &   9.98, 1.13      &   9.74, 0.73       &    0.14, 0.07    &   0.06, 0   &     0.09, 0                        \\ \hline

\multirow{3}{*}{$\sqrt{2}t_3$}    & \multicolumn{1}{l|}{$L_2$  error} &   2.88 (0.94)           &    2.89 (0.96)         &   2.67 (0.96)             &    0.48 (0.28)         &   0.51 (0.16)          &     0.54 (0.15)         \\
                           & \multicolumn{1}{l|}{$L_1$  error}&   14.64 (4.78)           &    14.87 (4.77)        &   13.74 (5)            &    1.54 (0.93)        &   1.61 (0.53)         &     1.72 (0.49)              \\
                           & \multicolumn{1}{l|}{MS}&   59.54 (11.57)           &    61.11 (11.62)        &   61.39 (9.38)            &    15.69 (1.13)        &   15.34 (1.17)         &     15.54 (3.05)              \\
                           & \multicolumn{1}{l|}{FPR, FNR}  &   9.22, 1.07         &   9.55, 1.47      &   9.59, 0.93       &    0.14, 0    &   0.07, 0   &     0.11, 0                        \\ \hline

\multirow{3}{*}{MixN}       & \multicolumn{1}{l|}{$L_2$  error} &   3.25 (0.87)           &    3.33 (0.94)         &   3.17 (0.93)             &    0.67 (0.35)         &   0.64 (0.22)          &     0.64 (0.2)         \\
                           & \multicolumn{1}{l|}{$L_1$  error}&   16.86 (4.64)           &    17.53 (5.17)        &   16.58 (5.01)            &    2.16 (1.3)        &   2.02 (0.73)         &     2.04 (0.69)              \\
                           & \multicolumn{1}{l|}{MS}&   61.23 (10.51)           &    62.36 (10.93)        &   62.55 (8.76)            &    15.87 (1.91)        &   15.29 (1.01)         &     15.24 (0.67)              \\
                           & \multicolumn{1}{l|}{FPR, FNR}  &   9.57, 1.27         &   9.82, 1.67      &   9.85, 1.33       &    0.18, 0    &   0.06, 0   &     0.05, 0                        \\ \hline

\multirow{3}{*}{LogNormal}   & \multicolumn{1}{l|}{$L_2$  error} &   3.68 (1.05)           &    3.64 (1)         &   3.4 (1.05)             &    0.9 (1.03)         &   0.64 (0.36)          &     0.74 (0.77)         \\
                           & \multicolumn{1}{l|}{$L_1$  error}&   18.76 (5.16)           &    19.08 (4.87)        &   17.72 (5.37)            &    2.95 (3.68)        &   2.03 (1.19)         &     2.43 (3.05)              \\
                           & \multicolumn{1}{l|}{MS}&   58.63 (10.39)           &    63.13 (11.89)        &   62.99 (8.07)            &    15.62 (1.72)        &   15.26 (0.69)         &     15.18 (0.66)              \\
                           & \multicolumn{1}{l|}{FPR, FNR}  &   9.12, 4.07         &   10.04, 3.73      &   9.98, 2.67       &    0.19, 1.93    &   0.06, 0.2   &     0.07, 1.07                        \\ \hline

\multirow{3}{*}{Weibull}    & \multicolumn{1}{l|}{$L_2$  error} &   3.01 (1.19)           &    2.83 (1.09)         &   2.66 (1.11)             &    0.75 (0.9)         &   0.59 (0.52)          &     0.64 (0.67)         \\
                           & \multicolumn{1}{l|}{$L_1$  error}&   15.09 (6.07)           &    14.89 (5.68)        &   13.66 (5.63)            &    2.5 (3.48)        &   1.94 (2.1)         &     2.11 (2.67)              \\
                           & \multicolumn{1}{l|}{MS}&   57.28 (11.85)           &    65.07 (9.4)        &   61.77 (7.9)            &    15.71 (1.21)        &   15.39 (0.91)         &     15.3 (0.98)              \\
                           & \multicolumn{1}{l|}{FPR, FNR}  &   8.78, 2.13         &   10.37, 1.53      &   9.69, 1.53       &    0.19, 1.27    &   0.09, 0.27   &     0.08, 0.67                        \\

\end{tabular}}}
\end{table}


In Example \ref{ex:3}, we only report results from the MCP-type PRAM estimators, since they have been shown to perform better than the Lasso-type estimators. In the homogeneous model with non-Gaussian covariates, both Table \ref{tb:contaminatex_inf} and Table \ref{tb:contaminatex_mallow} clearly indicate that the PRAM estimators with a well-chosen $w(\bx)$  perform better in all cases than those PRAM with $w(\bx)=1$. In addition, among those three weighted PRAM estimators, the weighted TA-MCP (WTA-MCP) and the weighted CA-MCP (WCA-MCP) again show advantages over the weighted HA-MCP (WHA-MCP) when the errors are heavy-tailed, which is consistent with the findings obtained in Example \ref{ex:1} and \ref{ex:2}. 

In conclusion, the PRAM estimator with a folded concave penalty (e.g. MCP penalty) renders promising performances in different settings, which is consistent with our theoretical results. Our simulation study also sheds some light on how to implement  robust high-dimensional M-estimators for real applications: when the data are strongly heavy-tailed or contaminated, regardless of asymmetry and/or heteroscedasticity, a re-descending PRAM estimator with a concave penalty yields better performances than a convex PRAM estimator in practice.


\begin{table}
  \caption{\label{tb:contaminatex_inf}Simulation results  under the homogeneous model with non-Gaussian covariates in Example \ref{ex:3}. The weighted PRAM estimators use $v(\bx) \equiv 1$ and $w(\bx)=\min\left\{1,\frac{4} {\|\bx\|_{\infty}}\right\}$. The mean $L_2$ error, $L_1$ error, MS, FPR (\%) and FNR (\%) out of 100 iterations are displayed. Standard deviations are listed in parentheses.}
 \centering
 \resizebox{\textwidth}{!}{%
 \fbox{%
\begin{tabular}{clrrrrrr}
                           & \multicolumn{1}{l|}{}                        & HA-MCP & WHA-MCP & TA-MCP & WTA-MCP  & CA-MCP & WCA-MCP \\ \hline
\multirow{3}{*}{N(0,4)}    & \multicolumn{1}{l|}{$L_2$  error} &   0.87 (0.91)           &    0.69 (0.61)         &   1.19 (1.58)             &    0.75 (0.83)         &   0.93 (1.1)          &     0.69 (0.6)         \\
                           & \multicolumn{1}{l|}{$L_1$  error}&   3.65 (4.29)           &    2.38 (2.71)        &   5.14 (7.6)            &    2.62 (3.77)        &   3.89 (5.17)         &     2.39 (2.7)              \\
                           & \multicolumn{1}{l|}{MS}&   36.92 (12.84)           &    17.7 (4.16)        &   36.88 (13.86)            &    17.69 (4.12)        &   36.07 (13.11)         &     17.89 (4.57)              \\
                           & \multicolumn{1}{l|}{FPR, FNR}  &   4.53, 0.2         &   0.56, 0.27      &   4.58, 2.27       &    0.58, 0.73    &   4.36, 0.67   &     0.6, 0.27                        \\ \hline

\multirow{3}{*}{$\sqrt{2}t_3$}    & \multicolumn{1}{l|}{$L_2$  error} &   0.91 (0.72)           &    0.65 (0.28)         &   1.11 (1.55)             &    0.63 (0.34)         &   0.87 (0.85)          &     0.61 (0.26)         \\
                           & \multicolumn{1}{l|}{$L_1$  error}&   3.75 (3.41)           &    2.12 (1.07)        &   4.8 (7.42)            &    2.07 (1.27)        &   3.56 (3.96)         &     1.98 (0.96)              \\
                           & \multicolumn{1}{l|}{MS}&   36.39 (11.54)           &    16.77 (2.78)        &   35.75 (11.68)            &    16.56 (2.63)        &   35.15 (11.86)         &     16.48 (2.71)              \\
                           & \multicolumn{1}{l|}{FPR, FNR}  &   4.42, 0.27         &   0.36, 0      &   4.36, 2.67       &    0.32, 0    &   4.16, 0.33   &     0.31, 0                        \\ \hline

\multirow{3}{*}{MixN}       & \multicolumn{1}{l|}{$L_2$  error} &   0.95 (0.89)           &    0.82 (0.71)         &   1.29 (1.52)             &    0.83 (0.75)         &   0.98 (0.99)          &     0.82 (0.74)         \\
                           & \multicolumn{1}{l|}{$L_1$  error}&   4.08 (4.25)           &    2.9 (3.31)        &   5.71 (7.33)            &    2.9 (3.4)        &   4.14 (4.69)         &     2.85 (3.4)              \\
                           & \multicolumn{1}{l|}{MS}&   38.22 (11.12)           &    18.5 (3.8)        &   39.42 (11.2)            &    17.9 (3.73)        &   37.65 (12.55)         &     17.88 (3.61)              \\
                           & \multicolumn{1}{l|}{FPR, FNR}  &   4.8, 0.47         &   0.74, 0.53      &   5.09, 1.8       &    0.61, 0.47    &   4.69, 0.8   &     0.61, 0.4                        \\ \hline

\multirow{3}{*}{LogNormal}   & \multicolumn{1}{l|}{$L_2$  error} &   2.26 (2.19)           &    1.31 (1.5)         &   2.74 (2.44)             &    1.38 (1.71)         &   2.02 (1.94)          &     1.36 (1.73)         \\
                           & \multicolumn{1}{l|}{$L_1$  error}&   10.24 (10.45)           &    4.8 (6.48)        &   12.51 (11.68)            &    5.19 (7.53)        &   9.07 (9.31)         &     4.85 (6.77)              \\
                           & \multicolumn{1}{l|}{MS}&   40.23 (11.37)           &    18.03 (4.55)        &   43.42 (12.09)            &    18.11 (4.35)        &   41.04 (12.25)         &     16.6 (3.6)              \\
                           & \multicolumn{1}{l|}{FPR, FNR}  &   5.39, 6         &   0.74, 3.67      &   6.16, 9.67       &    0.79, 4.93    &   5.52, 4.87   &     0.5, 5.6                        \\ \hline

\multirow{3}{*}{Weibull}    & \multicolumn{1}{l|}{$L_2$  error} &   1.6 (1.98)           &    1.11 (1.75)         &   1.84 (2.26)             &    0.92 (1.55)         &   1.44 (1.9)          &     0.85 (1.36)         \\
                           & \multicolumn{1}{l|}{$L_1$  error}&   7.11 (9.73)           &    4.34 (8.28)        &   8.04 (10.19)            &    3.34 (6.36)        &   6.17 (8.6)         &     3.03 (5.55)              \\
                           & \multicolumn{1}{l|}{MS}&   37.25 (11.37)           &    17.45 (5.57)        &   38.69 (13.33)            &    17.57 (4.88)        &   35.74 (10.84)         &     16.83 (3.61)              \\
                           & \multicolumn{1}{l|}{FPR, FNR}  &   4.7, 3.6         &   0.62, 3.8      &   5.08, 6.4       &    0.62, 2.93    &   4.38, 3.4   &     0.44, 2.13                        \\

\end{tabular}}
}
\end{table}

\begin{table}

  \caption{\label{tb:contaminatex_mallow}Simulation results under the homogeneous model with non-Gaussian covariates in Example \ref{ex:3}.  The weighted PRAM estimators use $v(\bx) \equiv 1$ and $w(\bx)=\min\left\{1,\frac{4} {\|\bx\|_{2}}\right\}$.  The mean $L_2$ error, $L_1$ error, MS, FPR (\%) and FNR (\%) out of 100 iterations are displayed. Standard deviations are listed in parentheses.}
  
 \centering
 \resizebox{\textwidth}{!}{%
 \fbox{%
\begin{tabular}{clrrrrrr}

                           & \multicolumn{1}{l|}{}                        & HA-MCP & WHA-MCP & TA-MCP & WTA-MCP  & CA-MCP & WCA-MCP \\ \hline
\multirow{3}{*}{N(0,4)}    & \multicolumn{1}{l|}{L2 loss} &   0.89 (1.04)           &    0.67 (0.5)         &   1.09 (1.43)             &    0.69 (0.59)         &   0.92 (1.1)          &     0.69 (0.59)         \\
                           & \multicolumn{1}{l|}{L1 loss}&   3.81 (5.09)           &    2.34 (2.65)        &   4.67 (6.71)            &    2.35 (2.69)        &   3.91 (5.43)         &     2.35 (2.69)              \\
                           & \multicolumn{1}{l|}{MSize}&   36.22 (12.68)           &    18.36 (4.73)        &   35.64 (11.15)            &    17.74 (3.87)        &   35.33 (12.3)         &     17.81 (4)              \\
                           & \multicolumn{1}{l|}{FP, FN}  &   4.39, 0.4         &   0.69, 0.07      &   4.31, 1.67       &    0.57, 0.2    &   4.21, 0.47   &     0.59, 0.2                        \\ \hline

\multirow{3}{*}{$\sqrt{2}t_3$}    & \multicolumn{1}{l|}{L2 loss} &   0.92 (0.77)           &    0.63 (0.21)         &   1.05 (1.42)             &    0.68 (0.81)         &   0.88 (0.88)          &     0.59 (0.17)         \\
                           & \multicolumn{1}{l|}{L1 loss}&   3.83 (3.66)           &    2.09 (0.84)        &   4.41 (6.38)            &    2.41 (4.53)        &   3.61 (3.98)         &     1.92 (0.63)              \\
                           & \multicolumn{1}{l|}{MSize}&   36.1 (12.25)           &    17.32 (3.23)        &   36.04 (11.49)            &    17.29 (4.23)        &   36.53 (13.48)         &     16.75 (2.91)              \\
                           & \multicolumn{1}{l|}{FP, FN}  &   4.36, 0.27         &   0.48, 0      &   4.41, 2.4       &    0.49, 0.53    &   4.45, 0.4   &     0.36, 0                        \\ \hline

\multirow{3}{*}{MixN}       & \multicolumn{1}{l|}{L2 loss} &   1.02 (0.98)           &    0.81 (0.77)         &   1.25 (1.47)             &    0.79 (0.73)         &   1 (1.07)          &     0.77 (0.74)         \\
                           & \multicolumn{1}{l|}{L1 loss}&   4.38 (4.69)           &    2.89 (3.55)        &   5.35 (6.62)            &    2.78 (3.43)        &   4.3 (5.14)         &     2.64 (3.19)              \\
                           & \multicolumn{1}{l|}{MSize}&   37.88 (11.43)           &    18.19 (3.79)        &   38.93 (11.55)            &    17.72 (4.52)        &   37.56 (11.68)         &     17.13 (2.62)              \\
                           & \multicolumn{1}{l|}{FP, FN}  &   4.74, 0.73         &   0.67, 0.53      &   5, 2.13       &    0.58, 0.53    &   4.68, 0.93   &     0.46, 0.67                        \\ \hline

\multirow{3}{*}{LogNormal}   & \multicolumn{1}{l|}{L2 loss} &   2.2 (2.08)           &    1.23 (1.45)         &   2.63 (2.47)             &    1.29 (1.58)         &   2.08 (2.03)          &     1.17 (1.43)         \\
                           & \multicolumn{1}{l|}{L1 loss}&   10.03 (10.28)           &    4.65 (6.64)        &   11.8 (11.43)            &    4.68 (6.4)        &   9.18 (9.52)         &     4.12 (5.52)              \\
                           & \multicolumn{1}{l|}{MSize}&   40.44 (12.02)           &    18.59 (5.94)        &   41.21 (12.11)            &    18.16 (4.19)        &   39.03 (12.73)         &     17.13 (3.91)              \\
                           & \multicolumn{1}{l|}{FP, FN}  &   5.42, 5.6         &   0.85, 3.47      &   5.72, 10.13       &    0.79, 4.33    &   5.11, 5.13   &     0.55, 3.67                        \\ \hline

\multirow{3}{*}{Weibull}    & \multicolumn{1}{l|}{L2 loss} &   1.59 (1.95)           &    1.16 (1.85)         &   1.71 (2.12)             &    0.97 (1.68)         &   1.49 (1.96)          &     0.93 (1.51)         \\
                           & \multicolumn{1}{l|}{L1 loss}&   7.18 (9.94)           &    4.65 (9.18)        &   7.46 (9.53)            &    3.81 (8.13)        &   6.53 (9.14)         &     3.39 (6.26)              \\
                           & \multicolumn{1}{l|}{MSize}&   37.86 (11.82)           &    17.59 (5.27)        &   38.83 (12.47)            &    18.09 (5.48)        &   38.12 (12.51)         &     17.18 (3.9)              \\
                           & \multicolumn{1}{l|}{FP, FN}  &   4.82, 3.4         &   0.65, 3.67      &   5.07, 5       &    0.72, 2.73    &   4.88, 3.8   &     0.53, 2.6                        \\ 

\end{tabular}}
}

\end{table}

\section{Data Example} \label{sec:realData}
In this section, we use the NCI-60 data, a gene expression data set collected from Affymetrix HG-U133A chip, to illustrate the performance of the PRAM estimators evaluated in Section {\ref{sec:simulation}}. The NCI-60 data consists of data from 60 human cancer cell lines and can be downloaded via the web application CellMiner (\url{http://discover.nci.nih.gov/cellminer/}). The study is to predict the protein expression on the KRT18 antibody from other gene expression levels. The expression level of the protein {\it keratin} 18 is known to be persistently expressed in carcinomas \cite{oshima1996oncogenic}. After removing the missing data, there are $n=59$ samples with  $21,944$ genes in the dataset. One can refer \cite{shankavaram2007transcript} for more details.

We perform some pre-screenings and keep only $p_1$
genes with largest variations and then choose $p_2$ genes out of them which are most correlated with the response variable. Here the final dataset is obtained by choosing $p_1=2000$ and $p_2=500$, yielding $n=59$ and $p=500$ for PRAM data analysis. Similar to our simulation studies, we then apply $6$ PRAM estimators to select important genes, with tuning parameters $\alpha$ and $\lambda$ chosen from the $10$-fold cross-validation. Since the TA-type and CA-type PRAM estimators perform similarly, we will only report results from four methods: HA-Lasso, CA-Lasso, HA-MCP and CA-MCP. 

\begin{figure}[!]
  	\centering{
  	\includegraphics[scale=0.75]{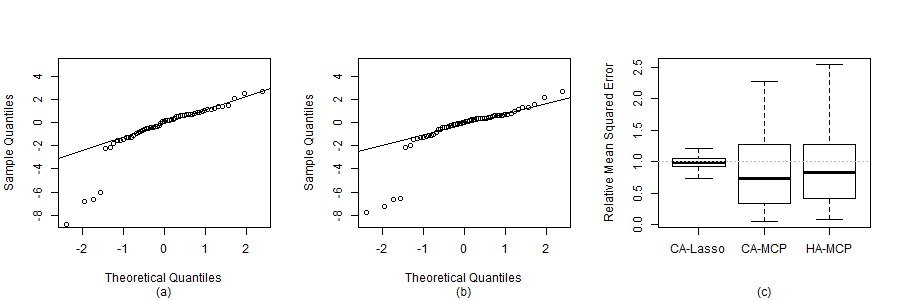}
  	 \caption{(a) The QQ plot of the residuals ($\hat{r}_i=y_i -  \bx_i^T\hat{\bb}$) from HA-MCP. (b) The QQ plot of the residuals from CA-MCP. (c) The boxplot of the relative mean squared prediction errors.}
  	 \label{fig:resid-box}
  	 }
 \end{figure}

The numbers of selected genes from four PRAM methods are 27 (HA-Lasso), 31 (CA-Lasso), 12 (HA-MCP), 5 (CA-MCP), respectively. HA-Lasso and CA-Lasso that selected 27 and 31 genes respectively could potentially result in over selection since the total sample size is only 59. Fig. \ref{fig:resid-box}(a) and Fig.  \ref{fig:resid-box}(b) show that the residual distributions generated from HA-MCP and CA-MCP both had a longer tail on the left side. It indicates that PRAM estimators with non-convex penalties can be resistant to the data contamination or data's irregularity due to the flexible robustness and nice variable selection property. 


For the sake of simplicity, we only report those selected genes and corresponding coefficient estimation by HA-MCP and CA-MCP in Table \ref{tb:vs}. According to our analysis, genes KRT8, NRN1 and GPX3 are selected by all four methods. It is not surprising for gene KRT8 since it has the largest correlation with the response variable and has a long history of being  paired with KRT18 in  cancer studies for cell death and survival, cellular growth and proliferation,  organismal injury and abnormalities, and so on \citep{li2016silac, Walker2007cancer}.
Gene NRN1 was investigated to be involved in melanoma migration, attachment independent growth, and vascular mimicry \citep{bosserhoff2017neurotrophin}. Recent studies showed that gene GPX3 plays as a tumor suppressor in lung cancer cell line \citep{an2018gpx3} and its down-regulation is related to pathogenesis of melanoma \citep{chen2016hypermethylation}.
Notice that gene ATP2A3 is also singled out by both HA-MCP and CA-MCP. This gene encodes the enzyme involved in calcium sequestration associated with muscular excitation and contraction, and was shown to act an important role in resveratrol anticancer activity in breast cancer cells \citep{izquierdo2017atp2a3}. In addition, Table \ref{tb:vs} indicates that gene GPNMB is only selected by CA-MCP. The GPNMB expression was found to be associated with a reduction in disease-free and overall survival in breast cancer and its over-expression had been identified in  numerous cancers \citep{maric2013glycoprotein}. Therefore, both genes (ATP2A3 and GPNMB) deserve further study in genetics research.


To further evaluate the prediction performance of those PRAM estimators, we randomly choose 6 observations as the test set and applied four methods to the rest patients to get the coefficients estimation, then compute the prediction error on the test set. We repeat the random splitting 100 times and the boxplots of the relative mean squared prediction error (RPE) with respect to HA-Lasso are shown in Fig. \ref{fig:resid-box}(c). A method with $\text{RPE}<1$ indicates a better performance than HA-Lasso. It is clearly seen from Fig. \ref{fig:resid-box}(c) that the MCP-type PRAM estimators have better predictions than those from the Lasso-type estimators, even though they select a much smaller number of variables. In addition, Fig. \ref{fig:resid-box}(c) together with Table \ref{tb:vs} show that a re-descending PRAM estimator with a non-convex penalty (e.g., CA-MCP) is more likely to give a more parsimonious model with better prediction performance, which is consistent with the findings from our simulation studies.

\begin{table}[]
  \caption{\label{tb:vs}Selected genes and the corresponding coefficient estimation by HA-MCP and CA-MCP. Probe IDs are listed in parentheses.}
 \resizebox{\textwidth}{!}{
\begin{tabular}{lllllll}
\hline\\
 HA-MCP      & KRT8            & NRN1            & GPX3         & CELF2           & CELF2           & LEF1            \\
       & (209008\_x\_at) & (218625\_at)    & (201348\_at) & (202156\_s\_at) & (202157\_s\_at) & (221558\_s\_at) \\
       & 6.230           & -1.505          & 0.031        & -0.002          & 0.000           & -0.003          \\
       & MEST            & FAR2            & PBX1         & CLEC11A         & CLEC11A         & ATP2A3          \\
       & (202016\_at)    & (220615\_s\_at) & (212148\_at) & (205131\_x\_at) & (211709\_s\_at) & (213036\_x\_at) \\
       & 0.009           & -0.037          & 0.035        & -0.036          & -0.017          & -0.003          \\
CA-MCP & KRT8            & NRN1            & GPX3         & GPNMB           & ATP2A3          &                 \\
       & (209008\_x\_at) & (218625\_at)    & (201348\_at) & (201141\_at)    & (213036\_x\_at) &                 \\
       & 6.122           & -0.775          & 0.693        & -0.556          & -0.763          &                 \\ \\\hline
\end{tabular}}
\end{table}

\section{Discussion}
The irregular settings including data asymmetry, heteroscedasticity and data contamination often exist due to the data high-dimensionality. It is very important to address these irregular settings both theoretically and numerically in high-dimensional data analysis. In this paper, we have proposed a class of  PRAM estimators for robust high-dimensional mean regression. The key feature of the PRAM is using a robust loss function (either a convex or a nonconvex redescending loss) with an adaptive robustification parameter to approximate the mean function produced from the quadratic loss. This approximation process can reduce the bias generated by data's irregularity in high-dimensional mean regression. The proposed framework is very general and it covers a wide range of loss functions and penalty functions, allowing both functions to be non-convex.

Theoretically, we establish statistical properties of PRAM in ultra high-dimensional settings when $p$ grows with $n$ at an almost exponential rate. Under sufficient conditions, we show its local estimation consistency at the minimax rate enjoyed by the LS-Lasso and further establish the local oracle properties of the PRAM estimators, including both selection consistency and asymptotic normality, when an amenable penalty is used. These theoretical results are applicable for general irregular settings, including the data are contaminated by heavy-tailed distribution and/or outliers in the random errors and covariates, and the random errors lack symmetry and/or homogeneity.  

To establish the estimation consistency and the oracle properties, the robustification parameter $\alpha$ should be chosen adaptively according to the sample size ($n$), dimension ($p$) and moment ($k$), as shown in Theorem \ref{Tm: est} and Theorem \ref{Tm: oracle}. In the presence of asymmetric and heavy-tailed/contaminated data, a well-behaved PRAM estimator can reduce the bias sufficiently (when $\alpha$ is large enough) and still maintain certain robustness (when $\alpha$ is not too large).

We also demonstrate a phase transition of conditional moment $k$ of the random error for PRAM estimators in ultra high-dimensional settings. Particularly, PRAM estimators have estimation consistency at a slower rate if $1<k<2$, but with a minimax rate enjoyed by the LS-Lasso if $k \ge 2$. For the oracle property, PRAM estimators can select variables consistently if only $k>1$. However, a larger $k$ permits the number of nonzero parameters to grow faster with $n$ and allows weaker minimum signal strength. Furthermore, the PRAM estimators obtain asymptotic normality if $k>2$.

Additionally, our numerical studies show satisfactory finite sample performances of PRAM estimators under irregular settings, which is consistent with our theoretical findings. In particular, our numerical results demonstrate the advantage of a re-descending PRAM estimator with a concave penalty such as TA-MCP and CA-MCP, using the HA-Lasso as the initial estimator, if the data are strongly heavy-tailed or contaminated.


Our research in this paper provides a systematic study of penalized M-estimation in high-dimensional mean regression. We hope this study shed some light on future directions of research, including devising similar theoretical guarantees for estimators with grouping structures in the covariates, or study of high-dimensional models with varying coefficients under general irregular settings.

\section*{Appendix}
\subsection*{A. Establishing the uniform RSC condition}
Let $\varepsilon_T= E\left[ P\left(|\epsilon_i| \ge \frac{T}{2}|\bx \right)\right]$ be the expected tail probability.
Below we establish some sufficient conditions where an unweighted $\mL_{\alpha, n}$ ($w(\bx) \equiv v(\bx)\equiv 1$) satisfies the uniform RSC condition in Assumptions \ref{as:RSC} with high probability. The uniform RSC condition for weighted loss can be established accordingly. 

\begin{theorem} \label{Tm: RSC}
Suppose $\mL_{\alpha,n}$ satisfies Assumption \ref{as:loss} and the covariate $\bx$ satisfies Assumption \ref{as:x}. If $ n \ge C_{10}s\log p$, then with probability at least $1-C_{11}\exp(-C_{12}\log p)$, the loss function $\mL_{\alpha,n}$ satisfies the Uniform RSC condition in  Assumption \ref{as:RSC} with 
\begin{equation*}
    \gamma=\frac{k_l}{32}, \quad  \quad \tau=\frac{C_{13}(3+ 2k_2)^2k_0^2T_0^2}{2r^2} \quad \text{and} \quad \alpha_0=\max\{ (2d_1)^{\frac{1}{k}},1\}\cdot T_0,
\end{equation*}
where $T_0 >0$ is a sufficiently large constant that satisfies 
\begin{equation} \label{ineq:RSC}
    C_{14}k_0^2\left(\sqrt{\varepsilon_{T_0}}+\exp\left(-\frac{C_{15}T_0^2}{k_0^2r^2}\right) \right) < \frac{k_l}{2+4k_2}.
\end{equation}

\end{theorem}

Theorem \ref{Tm: RSC} guarantees that the loss function $\mL_{\alpha, n}$ satisfies the uniform RSC condition with probability converging to 1. Note that the left-hand side of inequality (\ref{ineq:RSC}) is monotonically decreasing on $T_0$, meaning that inequality (\ref{ineq:RSC}) is always satisfied for a sufficiently large $T_0$. In addition, while keeping inequality (\ref{ineq:RSC}) satisfied, a larger $T_0$ (thus larger $\alpha_0$) actually allows a larger radius $r$ of local ball around $\bb^*$ and a more contaminated distribution of $\epsilon$. Theorem $\ref{Tm: RSC}$ implies that the Huber loss, Hampel loss, Tukey's biweight loss and Cauchy loss satisfy Assumption \ref{as:RSC} with high probability.

\subsection*{B. Proofs}
\label{proof}

{\bf Proof of Theorem \ref{Tm:1}}\\
	Let $l(x)=\frac{1}{2}x^2$. Observe that  
	\begin{equation*}
	\begin{split}
	E[\nabla \frac{w(\bx)}{v(\bx)}l((y-\bx^T\bb^*)v(\bx))]&=E[w(\bx)v(\bx)(y-\bx^T\bb^*)(-\bx)]\\
	&=E[w(\bx)v(\bx)\epsilon(-\bx)]\\
	&=E[E[\epsilon|\bx]w(\bx)v(\bx)(-\bx)]\\
	&=\b0,
	\end{split}
	\end{equation*}
	where the last equality follows from $E[\epsilon|\bx]=0$. Hence $\bb^*$ is the minimizer of $E[\frac{w(\bx)}{v(\bx)}l((y-\bx^T\bb)v(\bx))]$. Then it follows from Assumption \ref{as:x}(iii) that 
	\begin{equation} \label{eq:bound1}
	\begin{split}
		&  E[\frac{w(\bx)}{v(\bx)}l((y-\bx^T\bb_{\alpha}^*)v(\bx)) -\frac{w(\bx)}{v(\bx)}l((y-\bx^T\bb^*)v(\bx))]\\
		&=E\{w(\bx)v(\bx)[l(y-\bx^T\bb_{\alpha}^*)-l(y-\bx^T\bb^*)]\}\\
		&=\frac{1}{2}(\bb_{\alpha}^* - \bb^*)^TE[w(\bx)v(\bx)\bx\bx^T](\bb_{\alpha}^* - \bb^*) \ge \frac{1}{2}k_l \|\bb_{\alpha}^* - \bb^*\|_2^2
	\end{split}
	\end{equation}
	Let $g_{\alpha}(x)=l(x)- l_{\alpha}(x)$. Since $\bb_{\alpha}^*$ is the minimizer of $E [\frac{w(\bx)}{v(\bx)}l_{\alpha}((y-\bx^T\bb)v(\bx))]$ within a neighbour of $\bb^*$, we have
	\begin{equation} \label{ineq: exp_l_upper}
	\begin{split}
	&	E[\frac{w(\bx)}{v(\bx)}l((y-\bx^T\bb_{\alpha}^*)v(\bx)) - \frac{w(\bx)}{v(\bx)}l((y-\bx^T\bb^*)v(\bx))] \\
	=& E\{\frac{w(\bx)}{v(\bx)}[l((y-\bx^T\bb_{\alpha}^*)v(\bx))-l_{\alpha}((y-\bx^T\bb_{\alpha}^*)v(\bx))]\} + \\
	&E\{\frac{w(\bx)}{v(\bx)}[l_{\alpha}((y-\bx^T\bb_{\alpha}^*)v(\bx))- l_{\alpha}((y-\bx^T\bb^*)v(\bx))]\}+\\
	& E\{\frac{w(\bx)}{v(\bx)}[l_{\alpha}((y-\bx^T\bb^*)v(\bx))-l((y-\bx^T\bb^*)v(\bx))]\\
	\le & E[\frac{w(\bx)}{v(\bx)}g_{\alpha}((y-\bx^T\bb_{\alpha}^*)v(\bx))]-E[\frac{w(\bx)}{v(\bx)}g_{\alpha}((y-\bx^T\bb^*)v(\bx))]
	\end{split}
	\end{equation}

It follows from mean value theorem that
\begin{equation} \label{ineq:ga_diff}
\begin{split}
&E[\frac{w(\bx)}{v(\bx)}g_{\alpha}((y-\bx^T\bb_{\alpha}^*)v(\bx))-\frac{w(\bx)}{v(\bx)}g_{\alpha}((y-\bx^T\bb^*)v(\bx))]\\
=& E[w(\bx)\bx^T(\bb_{\alpha}^*-\bb^*) (z-l_{\alpha}'(z))] \\
\le & E[|w(\bx)\bx^T(\bb_{\alpha}^*-\bb^*)|| z-l_{\alpha}'(z)|]
\end{split}
\end{equation}
where $z=(y-\bx^T\tilde{\bb})v(\bx)$ and $\tilde{\bb}$ is a vector lying between $\bb^*$ and $\bb_{\alpha}^*$. Notice $l'_{\alpha}(0)=0$ in Assumption \ref{as:loss}(ii). By taking integral on each side of inequality  in Assumption \ref{as:loss}(iii), we have
\begin{equation} \label{ineq:l1bound}
    |u-l_{\alpha}'(u)| \le \frac{d_1}{k+1}|u|^{k+1}\alpha^{-k},
\end{equation}
for all $|u| \le \alpha$. Observe that
\begin{equation} \label{eq:I1+I2}
    \begin{split}
        E[|z-l_{\alpha}'(z)||\bx] =& E[|z-l_{\alpha}'(z)| \mathbbm{1}(|z| \le \alpha)|\bx] + E[|z-l_{\alpha}'(z)| \mathbbm{1}(|z| > \alpha)|\bx]\\
        =& I_1 + I_2.
    \end{split}
\end{equation}
From (\ref{ineq:l1bound}) we have
\begin{equation} \label{eq:I1}
\begin{split}
     I_1 =& E[|z-l_{\alpha}'(z)| \mathbbm{1}(|z| \le \alpha)|\bx]\\
     \le  &  \frac{d_1\alpha^{-k}}{k+1} E[ |z|^{k+1}\mathbbm{1}(|z| \le \alpha)|\bx] \\
     \le & \frac{d_1\alpha^{-k}}{k+1} E[ \frac{\alpha}{|z|}|z|^{k+1}|\bx] \\
     = & \frac{d_1\alpha^{1-k}}{k+1} E[|z|^{k}|\bx].
\end{split}
\end{equation}
Also observe that
\begin{equation} \label{eq:I2}
\begin{split}
    I_2 = & E[|z-l_{\alpha}'(z)| \mathbbm{1}(|z| > \alpha)|\bx]\\
    \le & E[|z| \mathbbm{1}(|z| > \alpha)|\bx] + E[|l_{\alpha}'(z)| \mathbbm{1}(|z| > \alpha)|\bx] \\
    < & \frac{1}{\alpha^{k-1}}E[|z|^{k}|\bx] + k_1\alpha E[\mathbbm{1}(|z| > \alpha)|\bx] \\
    = & \alpha^{1-k}E[|z|^{k}|\bx] + k_1 \alpha^{1-k}E[|z|^{k}|\bx] \\
    = & (1+k_1)\alpha^{1-k}E[|z|^{k}|\bx],
\end{split}
\end{equation}
where the second inequality follows from Assumption \ref{as:loss}(i).
Combining (\ref{eq:I1+I2}), (\ref{eq:I1}) and (\ref{eq:I2}), we obtain
\begin{equation} \label{eq: z-l'(z)1}
    E[|z-l_{\alpha}'(z)||\bx] \le (\frac{d_1}{k+1} +1 + k_1)\alpha^{1-k}E[|z|^{k}|\bx] =  C_1 \alpha^{1-k}E[|z|^{k}|\bx]
\end{equation}
where $C_1 = \frac{d_1}{k+1} +1 + k_1$ and $k$ is the constant that stated in Assumption \ref{as:loss}(iii), Assumption \ref{as:x}(i) and 3(ii).


Combining inequalities (\ref{ineq: exp_l_upper}), (\ref{ineq:ga_diff}) and (\ref{eq: z-l'(z)1}), we obtain
\begin{equation} \label{eq: bound2}
\begin{split}
&	E[\frac{w(\bx)}{v(\bx)}l((y-\bx^T\bb_{\alpha}^*)v(\bx)) - \frac{w(\bx)}{v(\bx)}l((y-\bx^T\bb^*)v(\bx))]\\\le& C_1\alpha^{1-k}E\{|y-\bx^T\tilde{\bb}|^{k}v(\bx)^k|w(\bx)\bx^T(\bb_{\alpha}^*-\bb^*)|\}\\
=&C_1\alpha^{1-k}E\{|\epsilon+\bx^T(\bb^*-\tilde{\bb})|^{k}v(\bx)^k|w(\bx)\bx^T(\bb_{\alpha}^*-\bb^*)|\}\\
\le & C_1(2/\alpha)^{k-1}\{ E[|\epsilon|^{k}v(\bx)^k|w(\bx)\bx^T(\bb_{\alpha}^*-\bb^*)|] +\\ &E[|\bx^T(\bb^*-\tilde{\bb})|^{k}v(\bx)^k|w(\bx)\bx^T(\bb_{\alpha}^*-\bb^*)|]\},
\end{split}
\end{equation}
where the last inequality follows from Minkowski inequality. Note that
\begin{equation} \label{ineq:left}
\begin{split}
E[|\epsilon|^{k}v(\bx)^k|w(\bx)\bx^T(\bb_{\alpha}^*-\bb^*)|]=&E[E(|\epsilon|^{k}|\bx)v(\bx)^k|w(\bx)\bx^T(\bb_{\alpha}^*-\bb^*)|]\\
\le & \{E[E(|\epsilon|^{k}|\bx)v(\bx)^k]^2\}^{\frac{1}{2}} \{E[w(\bx)\bx^T(\bb_{\alpha}^*-\bb^*)]^2\}^{\frac{1}{2}}\\
\le &\sqrt{M_kk_u}\|\bb_{\alpha}^*-\bb^*\|_2,
\end{split}
\end{equation}
where the first inequality follows from H\"{o}lder inequality and the last inequality follows from Assumption \ref{as:x}(i) and (iii). Observe that,
\begin{equation} \label{ineq:right}
\begin{split}
E[|\bx^T(\bb^*-\tilde{\bb})|^{k}v(\bx)^k|w(\bx)\bx^T(\bb_{\alpha}^*-\bb^*)|] \le &
\{E[v(\bx)\bx^T(\bb^*-\tilde{\bb})]^{2k}\}^{\frac{1}{2}} \{E[w(\bx)\bx^T(\bb_{\alpha}^*-\bb^*)]^2\}^{\frac{1}{2}}\\
\le & R_0^k\sqrt{q_{k}k_u}\|\bb_{\alpha}^*-\bb^*\|_2,
\end{split}
\end{equation}
where $R_0$ is defined in (\ref{eq:beta_alpha_star}) and  the last inequality follows from Assumption \ref{as:x}(ii) and \ref{as:x}(iii). By inequalities (\ref{eq:bound1}), (\ref{eq: bound2}),  (\ref{ineq:left}), (\ref{ineq:right}) we have
\begin{equation*}
\|\bb_{\alpha}^*-\bb^*\|_2 \le 2^k C_1k_l^{-1}\sqrt{k_u}(\sqrt{M_k}+R_0^k\sqrt{q_k}){\alpha}^{1-k}.
\end{equation*}
\QEDB\\

To prove Theorem \ref{Tm: est}, we need the following Lemma \ref{Lemma: gradient_bound}.

\begin{lemma} \label{Lemma: gradient_bound}
Suppose $\mL_{\alpha, n}$ in (\ref{eq:GPRAM-estimator}) satisfies Assumption \ref{as:loss} and the random errors and covariates satisfy Assumption \ref{as:x}. For any $t \in (0,n)$, let $\max\{(\frac{2d}{R_0})^{\frac{1}{k-1}},C_2 (\frac{n}{t})^{\frac{\zeta}{k-1}}\} \le \alpha \le C_3(\frac{n}{t})^{1-\zeta}$ with $\zeta \le \min\{\frac{k-1}{k}, \frac{1}{2}\}$. Then we have
\begin{equation*}
    \|\nabla\mL_{\alpha, n}(\bb^*)\|_\infty \le C_{50}(\frac{t}{n})^{\zeta}
\end{equation*}
with probability at least $1-2p\exp(-t)$.
\end{lemma}
\textit{Proof.} 
	The gradient of $\mL_{\alpha,n}$ is
	\begin{equation}
	\begin{split}
	\nabla \mL_{\alpha,n}(\bb^*)=&-\frac{1}{n} \sum_{i=1}^{n} w(\bx_i)l_{\alpha}'((y_i-\bx_i^T\bb^*)v(\bx_i))\bx_i.
	\end{split}
	\end{equation}
Recall $\bb_{\alpha}^*$ is the minimizer of $E [\frac{w(\bx)}{v(\bx)}l_{\alpha}((y-\bx^T\bb)v(\bx))]$ within a neighbour of $\bb^*$ defined in (\ref{eq:beta_alpha_star}). When $\alpha \ge (\frac{2d}{R_0})^\frac{1}{k-1}$ where $d=2^k C_1k_l^{-1}\sqrt{k_u}(\sqrt{M_k}+R_0^k\sqrt{q_k})$, we have $\|\bb_{\alpha}^*-\bb^*\|_2 \le \frac{R_0}{2} < R_0$ under the result of Theorem \ref{Tm:1}. Hence $\bb_{\alpha}^*$ is an interior point of program (\ref{eq:beta_alpha_star}). Then we have $E[w(\bx)l_{\alpha}'((y-\bx^T\bb_{\alpha}^*)v(\bx))\bx]=\b0$. Observe that
\begin{equation} \label{eq:meanBound}
\begin{split}
E[w(\bx_i)l_{\alpha}'((y_i-\bx_i^T\bb^*)v(\bx_i))x_{ij}] = & E[w(\bx_i)l_{\alpha}'((y_i-\bx_i^T\bb^*)v(\bx_i))x_{ij}] - E[w(\bx_i)l_{\alpha}'((y_i-\bx_i^T\bb_{\alpha}^*)v(\bx))x_{ij}]\\
 \le & k_2E[|v(\bx_i)\bx_i^T(\bb_{\alpha}^*-\bb^*)||w(\bx_i)x_{ij}|]\\
\le & k_2\{E|v(\bx_i)\bx_i^T(\bb_{\alpha}^*-\bb^*)|^2\}^{\frac{1}{2}}\{E|w(\bx_i)x_{ij}|^2\}^{\frac{1}{2}}\\
\le & k_2\sqrt{q_1}\|\bb_{\alpha}^*-\bb^*\|_2\sqrt{k_0^2+d_2^2}\\
\le & d_3 \alpha^{1-k},
\end{split}
\end{equation} 
where $\max_{1 \le j \le p}|E[w(\bx_i)\bx_{ij}]| < d_2 <\infty$ and  $d_3=2^kk_2\sqrt{q_1(k_0^2+d_2^2)k_u}C_1k_l^{-1}(\sqrt{M_k}+2^kR_0^k\sqrt{q_k})$. Note that the first inequality is from Assumption \ref{as:loss}(ii) and the third inequality follows from Assumption \ref{as:x}(ii) and (iv). And the last inequality is from Theorem \ref{Tm:1}. 

Next we  bound the $E|w(\bx_i)l_{\alpha}'((y_i-\bx_i^T\bb^*)v(\bx_i))x_{ij}|^m$ from the above, for  $m \ge 2$. Define $\tilde{k}=\min\{k,2\}$.
By Assumption \ref{as:loss} and \ref{as:x}(i)   we have
\begin{equation} \label{ineq:el1m}
\begin{split}
E|w(\bx_i)l_{\alpha}'(\epsilon_iv(\bx_i))x_{ij}|^m \le& E[(k_1\alpha)^{m-\tilde{k}}(k_2\epsilon_iv(\bx_i))^{\tilde{k}}|w(\bx_i)x_{ij}|^m]\\
\le & k_1^{m-\tilde{k}}\alpha^{m-\tilde{k}}k_2^{\tilde{k}}E[(\epsilon_iv(\bx_i))^{\tilde{k}}|w(\bx_i)x_{ij}|^m]\\
\le & k_1^{m-\tilde{k}}\alpha^{m-\tilde{k}}k_2^{\tilde{k}}\{ E[E(\epsilon_i^{\tilde{k}}|\bx_i)v(\bx_i)^{\tilde{k}}]^2\}^{1/2}\{E[(w(\bx_i)x_{ij})^{m}]^2 \}^{1/2}\\
\le & k_1^{m-\tilde{k}}\alpha^{m-\tilde{k}}k_2^{\tilde{k}}\sqrt{M_{\tilde{k}}} \{E[(w(\bx_i)x_{ij})^{m}]^2 \}^{1/2}.
\end{split}
\end{equation}

Let $\mu_j = E[w(\bx_i)\bx_{ij}]$, $j=1,2,\dots, p$. Then we have
\begin{equation} \label{eq:boundwxm}
    \begin{split}
        E|w(\bx_i)\bx_{ij}|^m=& E|w(\bx_i)\bx_{ij}-\mu_j+\mu_j|^m\\
        \le & E[2^{m-1}(|w(\bx_i)\bx_{ij}-\mu_j|^m+|\mu_j|^m)]\\
        \le & 2^{m-1} [E|w(\bx_i)\bx_{ij}-\mu_j|^m + d_2^m]\\
        \le & 2^{m-1} [m(\sqrt{2})^{m}k_0^m \Gamma(\frac{m}{2}) + d_2^m],
    \end{split}
\end{equation}
where the last inequality follows from Assumption \ref{as:x}(iv), by which $w(\bx_i)x_{ij}$ is sub-Gaussian hence for $m > 0$(\cite{rivasplata2012subgaussian})
\begin{equation*}\label{eq:gaussian}
E|w(\bx_i)x_{ij}-\mu_j|^m \le m(\sqrt{2})^{m}k_0^m\Gamma(\frac{m}{2}).
\end{equation*}

By taking $m=2$ in (\ref{ineq:el1m}), we have
\begin{equation}  \label{ineq:psi_squred}
\begin{split}
E[w(\bx_i)l_{\alpha}'((y_i-\bx_i^T\bb^*)v(\bx_i))x_{ij}]^2 \le & k_1^{2-\tilde{k}}\alpha^{2-\tilde{k}}k_2^{\tilde{k}}\sqrt{M_{\tilde{k}}}\{E[(w(\bx_i)x_{ij})^{2}]^2 \}^{1/2}\\
\le &  k_1^{2-\tilde{k}}\alpha^{2-\tilde{k}}k_2^{\tilde{k}}\sqrt{M_{\tilde{k}}}(128k_0^4 + 8d_2^4)^\frac{1}{2}\\
\le & d_4 \alpha^{2-\tilde{k}},
\end{split}
\end{equation}
where $d_4=\sqrt{2}k_1^{2-\tilde{k}}k_2^{\tilde{k}}\sqrt{M_{\tilde{k}}}(8k_0^2+2d_2^2) $ and the second inequality follows from (\ref{eq:boundwxm}).

For $m \ge 3$, by replacing $m$ by $2m$ in (\ref{eq:boundwxm}), we obtain
\begin{equation} \label{ineq: wx2mbound}
    \begin{split}
     \{E|w(\bx_i)\bx_{ij}|^{2m} \}^{\frac{1}{2}} \le & \{2^{2m-1} (2m)2^{m}k_0^{2m} \Gamma(m) + 2^{2m-1} d_2^{2m}\}^{\frac{1}{2}}\\
         \le& 2^{\frac{3m}{2}} k_0^m\sqrt{m!} + 2^{m-\frac{1}{2}}d_2^m\\
         =& (2^{\frac{3m}{2}}k_0^m \frac{2}{\sqrt{m!}} + \frac{2^{m+\frac{1}{2}}d_2^m}{m!})\frac{m!}{2}\\
         \le & (2^{\frac{3m}{2}}k_0^m + 2^{m-1}d_2^m)\frac{m!}{2}\\
         = & [(2^{\frac{3}{2}}k_0)^{m-2} \cdot (2^{\frac{3}{2}}k_0)^2 + (2d_2)^{m-2} \cdot 2d_2^2]\frac{m!}{2} \\
         \le& \frac{m!}{2}(2^{\frac{3}{2}}k_0 + 2d_2)^{m-2}(8k_0^2 + 2d_2^2).
    \end{split}
\end{equation}

Combining inequality (\ref{ineq:el1m}) and (\ref{ineq: wx2mbound}), we have
\begin{equation*}
\begin{split}
E|w(\bx_i)l_{\alpha}'(\epsilon_iv(\bx_i))x_{ij}|^m \le &  k_1^{m-\tilde{k}}\alpha^{m-\tilde{k}}k_2^{\tilde{k}}\sqrt{M_{\tilde{k}}}[\frac{m!}{2}(2^{\frac{3}{2}}k_0 + 2d_2)^{m-2}(8k_0^2 + 2d_2^2)]\\
< &\frac{m!}{2}(4(k_0+d_2)k_1\alpha)^{m-2}(k_1^{2-\tilde{k}}\alpha^{2-\tilde{k}}k_2^{\tilde{k}}\sqrt{M_{\tilde{k}}}(8k_0^2+2d_2^2))\\
< &\frac{m!}{2}(4(k_0+d_2)k_1\alpha)^{m-2}d_4\alpha^{2-\tilde{k}},
\end{split}
\end{equation*}
By Bernstein inequality (Proposition 2.9 of \cite{massart2007concentration}) we have
\begin{equation*}
\begin{array}{ll}
&P\left(|\frac{1}{n}\sum_{i=1}^{n}w(\bx_i)l_{\alpha}'((y_i-\bx_i^T\bb^*)v(\bx_i))x_{ij}-\frac{1}{n}\sum_{i=1}^{n}E[w(\bx_i)l_{\alpha}'((y_i-\bx_i^T\bb^*)v(\bx_i))x_{ij}]|\right.\\
&\quad\quad\quad\ge \left. \sqrt{\frac{2d_4t\alpha^{2-\tilde{k}}}{n}}+\frac{4(k_0+d_2)k_1 t\alpha}{n} \right)\\
 &\le 2\exp(-t),
\end{array}
\end{equation*}
for any $t>0$. It implies that
\begin{equation*}
\begin{array}{ll}
&P\left((|\frac{1}{n}\sum_{i=1}^{n}w(\bx_i)l_{\alpha}'((y_i-\bx_i^T\bb^*)v(\bx_i))x_{ij}|\right.\\
&\quad\quad\quad\ge  \left.\sqrt{\frac{2d_4t\alpha^{2-\tilde{k}}}{n}}+\frac{4(k_0+d_2)k_1 t\alpha}{n}+|\frac{1}{n}\sum_{i=1}^{n}E[w(\bx_i)l_{\alpha}'((y_i-\bx_i^T\bb^*)v(\bx_i))x_{ij}]|\right)\\ 
&\le 2\exp(-t).
\end{array}
\end{equation*}
By the bound in (\ref{eq:meanBound}),
\begin{equation} \label{eq:concentration}
\begin{array}{ll}
P(|\frac{1}{n}\sum_{i=1}^{n}w(\bx_i)l_{\alpha}'((y_i-\bx_i^T\bb^*)v(\bx_i))x_{ij}|
\ge  \sqrt{\frac{2d_4t\alpha^{2-\tilde{k}}}{n}}+\frac{4(k_0+d_2)k_1 t\alpha}{n}+d_3\alpha^{1-k}) \le 2\exp(-t).
\end{array}
\end{equation}
Let $\zeta \le \min\{\frac{k-1}{k}, \frac{1}{2}\}$ and $t < n$. 
By conditions we have 
\begin{equation} \label{eq:alpha}
C_2(\frac{n}{t})^{\frac{\zeta}{k-1}} \le \alpha \le C_3 (\frac{n}{t})^{1-\zeta}.
\end{equation}
If $k \ge 2$, then
\begin{equation*}
    \sqrt{\frac{2d_4t\alpha^{2-\tilde{k}}}{n}} \le \sqrt{2d_4}(\frac{t}{n})^{\frac{1}{2}} \le \sqrt{2d_4} (\frac{t}{n})^{\zeta}.
\end{equation*}
If $1 < k <2$, then 
\begin{equation*} 
     \sqrt{\frac{2d_4t\alpha^{2-\tilde{k}}}{n}}=\sqrt{\frac{2d_4t\alpha^{2-k}}{n}} \le \sqrt{2d_4}C_3^{\frac{2-k}{2}}(\frac{t}{n})^{\frac{(2-k)\zeta-1+k}{2}} \le \sqrt{2d_4}C_3^{\frac{2-k}{2}} (\frac{t}{n})^{\zeta},
\end{equation*}
where the first inequality follows from the second inequality in (\ref{eq:alpha}). Hence, we can choose some $C_3>1$ such that
\begin{equation} \label{ineq:lambdan_1}
     \sqrt{\frac{2d_4t\alpha^{2-\tilde{k}}}{n}}\le \sqrt{2d_4}C_3 (\frac{t}{n})^{\zeta} 
\end{equation}
By equation (\ref{eq:alpha}) we have
\begin{equation}
    \frac{4(k_0+d_2)k_1 t\alpha}{n} \le 4C_3(k_0+d_2)k_1 (\frac{t}{n})^\zeta,
\end{equation}
and 
\begin{equation}
    d_3\alpha^{1-k} \le C_2d_3 (\frac{t}{n})^\zeta.
\end{equation}
Together with (\ref{ineq:lambdan_1}) we obtain
\begin{equation*}
\sqrt{\frac{2d_4t\alpha^{2-\tilde{k}}}{n}}+\frac{4(k_0+d_2)k_1 t\alpha}{n}+d_3\alpha^{1-k} \le C_{50}(\frac{t}{n})^\zeta,
\end{equation*}
where $C_{50}=\sqrt{2d_4}C_3 + 4C_3(k_0+d_2)k_1 +C_2d_3$.
Hence by (\ref{eq:concentration}), it gives
\begin{equation} \label{ineq: 1d1}
P\left (|\frac{1}{n}\sum_{i=1}^{n}w(\bx_i)l_{\alpha}'((y_i-\bx_i^T\bb^*)v(\bx_i))x_{ij}|
\ge C_{50}(\frac{t}{n})^\zeta \right ) \le 2\exp(-t).
\end{equation}
It then follows from union inequality that 
\begin{equation} \label{ineq:pd1}
P\left (\|\nabla \mL_{\alpha,n}(\bb^*)\|_{\infty}
\ge   C_{50}(\frac{t}{n})^\zeta \right) \le 2p\exp(-t).
\end{equation}

\QEDB\\

{\bf Proof of Theorem \ref{Tm: est}}\\

By Lemma \ref{Lemma: gradient_bound} with $t=(1+C_4)\log p$, we obtain
\begin{equation} \label{eq:gradient_bound}
    P\left (\|\nabla \mL_{\alpha,n}(\bb^*)\|_{\infty}
\ge   C_5(\frac{\log p}{n})^\zeta \right) \le 2\exp(-C_4\log p),
\end{equation}
where $\alpha$ satisfies
\begin{equation*}
\max\{(\frac{2d}{R_0})^{\frac{1}{k-1}},C_{21}(\frac{n}{\log p})^{\frac{\zeta}{k-1}}\} \le \alpha \le C_{31} (\frac{n}{\log p})^{1-\zeta},
\end{equation*}
$C_5=C_{50}(1+C_4)^\zeta$, $C_{21}=C_2(1+C_4)^{-\frac{\zeta}{k-1}}$ and $C_{31}=C_3(1+C_4)^{\zeta-1}$. And the rest of the result follows immediately from Theorem 1 in \cite{loh2017statistical}.
In particular, similar to the proof of Theorem 1 in \cite{loh2017statistical}, we first restrict the optimization within the uniform RSC region and defined $\hat{\bb} \in \RR^p$ such that
\begin{equation} \label{eq:rsc_program}
    \hat{\bb}\in \argmin_{\|\bb\|_1\le R, \; \|\bb-\bb^*\|_2\le r} \left\{ \mL_{\alpha,n}(\bb)+\rho_\lambda(\bb)\right\}.
\end{equation}
Then $\hat{\bb}$ is a stationary point of the program (\ref{eq:rsc_program}). By equation (\ref{eq:gradient_bound}) and following the same arguments of the proof of Theorem 1 in \cite{loh2017statistical}, we obtain
	\begin{equation*}
	   \|\hat{\bb}-\bb^*\|_2 \le \frac{24C_5}{4\gamma-3\mu}\sqrt{s} \left(\frac{\log p}{n}\right)^\zeta  \text{ and  } \|\hat{\bb}-\bb^*\|_1 \le \frac{96C_5}{4\gamma-3\mu}s\left (\frac{\log p}{n} \right)^\zeta,
	\end{equation*}
with probability $1-2\exp(-C_4\log p)$. Let $C_0=(\frac{24C_5}{4\gamma-3\mu})^{1/\zeta}$. Given $n > C_0r^{-\frac{1}{\zeta}}s^{\frac{1}{2\zeta}}\log p$, we have
\begin{equation*}
    \|\hat{\bb}-\bb^*\|_2 \le \frac{24C_5}{4\gamma-3\mu}\sqrt{s} \left(\frac{\log p}{n}\right)^\zeta < r,
\end{equation*}
which means that the point $\hat{\bb}$ will lie in the interior of the sphere of radius $r$ around $\bb^*$. Therefore,  $\hat{\bb}$ is also a stationary point of the original program (\ref{eq:estimator}) and the existence of such local stationary points is guaranteed.

\QEDB\\

To prove Theorem \ref{Tm: oracle}, we need the following result adopted directly from the Lemma 1 in \cite{loh2017statistical}. 
\begin{lemma} \label{lemma:restrict-cov}
Suppose $\mL_{\alpha, n}$ satisfies the local RSC condition (\ref{as:RSC}) and $n \ge \frac{2\tau}{\gamma} s \log p$. Then $\mL_{\alpha, n}$ is strongly convex over the region $S_r = \{ \bb \in \RR^p: supp(\bb) \subseteq S, \|\bb - \bb^*\|_2 \le r\}$.
\end{lemma}
Proof. The proof is similar to the proof of Lemma 1 in \cite{loh2017statistical}. 

\QEDB\\

{\bf Proof of Theorem \ref{Tm: oracle}}\\
The proof is an adaptation of the arguments of Theorem 2 in the paper  \cite{loh2017statistical}. We follow the three steps of the primal-dual witness (PDW) construction described in that paper:
\begin{enumerate}
    \item[(i)] Optimize the restricted program 
    \begin{equation} \label{eq:restircted}
        \hat{\bb}_S \in \argmin_{\bb \in \bb^S: \|\bb\|_1 \le R}  \left\{  \mL_{\alpha, n}(\bb) + \rho_{\lambda}(\bb)
        \right \},
    \end{equation}
    and establish that $\|\hat{\bb}_S\|_1 < R$.
    
    \item[(ii)] Recall $q_\lambda(\bb) = \lambda\|\bb\|_1-\rho_\lambda(\bb)$ defined in Section \ref{sec:4}. Define $\hat{\bz}_S \in \partial \|\hat{\bb}_S\|_1$, and choose $\hat{ \bz}= (\hat{\bz}_S, \hat{\bz}_{S^c})$ to satisfy the zero-subgradient condition
    \begin{equation} \label{eq:zero-subgradient}
        \nabla \mL_{\alpha, n}(\hat{\bb}) - \nabla q_{\lambda}(\hat{\bb})+\lambda \hat{\bz} = \b0,
    \end{equation}
    where $\hat{\bb} = (\hat{\bb}_S, \pmb 0_{S^c})$. Show that $\hat{\bb}_S = \hat{\bb}^\bo_S$ and establish strict dual feasibility: $\|\hat{\bz}_{S^c}\|_\infty < 1$.
    
    \item[(iii)] Verify via second order conditions that $\hat{\bb}$ is a local minimum of the program (\ref{eq:GPRAM-estimator}) and conclude that all stationary points $\tilde{\bb}$ satisfying $\|\tilde{\bb} - \bb^*\|_2 \le r$ are supported on S and agree with $\hat{\bb}^\bo$.
\end{enumerate}

\textbf{Proof of Step (i) }:  By applying Theorem \ref{Tm: est} to the restricted program (\ref{eq:restircted}), we have
\begin{equation*}
    \|\hat{\bb}_S-\bb^*_S\|_1 \le \frac{96\lambda s}{4\gamma-3\mu},
\end{equation*}
and thus
\begin{equation*}
    \|\hat{\bb}_S\|_1 \le \|\bb^*\|_1 + \|\hat{\bb}_S-\bb^*_S\|_1 \le \frac{R}{2} +  \frac{96\lambda s}{4\gamma-3\mu} < R,
\end{equation*}
under the assumption of the theorem. This complete step (i) of the PDW construction. 

\QEDB\\

To prove step (ii), we need the following Lemma \ref{Lemma: restrict-bound} and \ref{Lemma: restricted-strict-convex}:
\begin{lemma}  \label{Lemma: restrict-bound}
Under the conditions of Theorem \ref{Tm: oracle}, we have the bound 
\begin{equation*}
     \|\hat{\bb}^\bo_S-\bb^*_S\|_2 \le C_6s^{-1/2}(\frac{\log p}{n})^{\zeta}
\end{equation*}
and $\hat{\bb}_S = \hat{\bb}^\bo_S$ with probability at least $1-2 \exp(-C_{41}s^{-\frac{1}{\zeta}}\log p)$.
\end{lemma}
\textit{Proof.} 
Recall $\hat{\bb}^\bo = (\hat{\bb}^\bo_S, \pmb 0_{S^c})$. By the optimality of the oracle estimator in (\ref{eq:oracle}), we have
\begin{equation} \label{ineq:oracle}
    \mL_{\alpha, n} (\hat{\bb}^\bo) \le \mL_{\alpha, n}(\bb^*).
\end{equation}
When $n \ge \frac{2\tau}{\gamma}s\log p$, by Lemma \ref{lemma:restrict-cov}, $\mL_{\alpha, n}(\bb)$ is strongly convex over restricted region $S_r=\{\|\bb - \bb^*\|_2 \le r\}$ . Hence,
\begin{equation}
    \mL_{\alpha, n}(\bb^*) + \langle \nabla \mL_{\alpha, n}(\bb^*), \hat{\bb}^\bo - \bb^* \rangle + \frac{\gamma}{4} \|\hat{\bb}^\bo - \bb^*\|_2^2 \le \mL_{\alpha, n}(\hat{\bb}^\bo).
\end{equation}
Together with inequality (\ref{ineq:oracle}) we obtain
\begin{equation*}
\begin{array}{ll}
\frac{\gamma}{4}  \|\hat{\bb}^\bo - \bb^*\|_2^2 &\le \langle \nabla \mL_{\alpha, n}(\bb^*), \bb^* - \hat{\bb}^\bo \rangle \le \| \nabla (\mL_{\alpha, n}(\bb^*))_S\|_\infty \cdot \|\hat{\bb}^\bo - \bb^*\|_1\\
& \le \sqrt{s} \| \nabla (\mL_{\alpha, n}(\bb^*))_S\|_\infty \cdot \|\hat{\bb}^\bo - \bb^*\|_2, 
\end{array}
\end{equation*}
implying that 
\begin{equation} \label{ineq:res-l2-bound}
    \|\hat{\bb}^\bo - \bb^*\|_2 \le \frac{4\sqrt{s}}{\gamma}\| \nabla (\mL_{\alpha, n}(\bb^*))_S\|_\infty.
\end{equation}
By Lemma \ref{Lemma: gradient_bound} we have
\begin{equation*}
P\left (\|\nabla \mL_{\alpha,n}(\bb^*_{S})\|_{\infty}
\ge   C_{50}(\frac{t}{n})^\zeta \right) \le 2s\exp(-t).
\end{equation*}
with $C_2(\frac{n}{t})^{\frac{\zeta}{k-1}} \le \alpha \le C_3 (\frac{n}{t})^{1-\zeta}$. Let $t=C_{51}s^{-\frac{1}{\zeta}}\log p$, we obtain
\begin{equation} \label{ineq:res-tm1}
    \|(\nabla \mL_{\alpha, n}(\bb^*))_S\|_\infty = \|\nabla \mL_{\alpha, n}(\bb^*_S)\|_\infty \le C_{50}C_{51}^{\zeta} s^{-1}(\frac{\log p}{n})^\zeta,
\end{equation}
with probability at least $1-2 \exp(-C_{41}s^{-\frac{1}{\zeta}}\log p)$, where we require $s^{\frac{1}{\zeta}}\log s =\bo (\log p)$ and $\alpha$ satisfies
\begin{equation}\label{ineq: alpha_range2}
C_{22}s^{\frac{1}{k-1}}(\frac{n}{\log p})^{\frac{\zeta}{k-1}} \le \alpha \le C_{32}s^{\frac{1-\zeta}{\zeta}} (\frac{n}{\log p})^{1-\zeta}
\end{equation}
Combining inequality (\ref{ineq:res-l2-bound}) and (\ref{ineq:res-tm1}), we obtain 
\begin{equation} \label{ineq: res-inf-bound}
     \|\hat{\bb}^\bo - \bb^*\|_2 \le  C_6s^{-1/2}(\frac{\log p}{n})^{\zeta}
\end{equation}
as desired, where $C_6=4C_{50}C_{51}^{\zeta}/\gamma$. 

Next we show $\hat{\bb}_S = \hat{\bb}^\bo_S$. When $n > C_6^2r^{-2}s^{-1}\log p $, we have $\|\hat{\bb}^\bo_S - \bb^*_S\|_2 < r$ and thus $\hat{\bb}^\bo_S$ is an interior point of the oracle program in (\ref{eq:oracle}), implying
\begin{equation} \label{eq:oracle-gradient}
    \nabla \mL_{\alpha, n}(\hat{\bb}^\bo_S) =\b0.
\end{equation}
By assumption we have $\lambda=C_{71}(\frac{\log p}{n})^\zeta$ and
$\beta^*_{\min} \ge C_{72}(\frac{\log p}{n})^{\zeta}$, where we choose $C_{72}=C_6+C_{71}\delta$. Together with inequality (\ref{ineq: res-inf-bound}), we have
\begin{equation*}
    \begin{array}{ll}
    |\hat{\beta}^\bo_j| \ge |\beta^*_j| - |\hat{\beta}^\bo_j - \beta^*_j| &\ge \beta^*_{\min} - \|\hat{\bb}^\bo_S-\bb^*_S\|_\infty \\
    & \ge ( C_6+C_{71}\delta)(\frac{\log p}{n})^\zeta -C_6(\frac{ \log p}{n})^\zeta \\
    & = \delta \lambda.
    \end{array}
\end{equation*}
for all $j \in S$. Together with the assumption that $\rho_{\lambda}$ is $(\mu, \delta)$-amenable, that is, Assumption \ref{as:penalty}(vii), we have 
\begin{equation} \label{eq:oracle-penalty}
    \nabla q_{\lambda}(\hat{\bb}^\bo_S) = \lambda \text{sign}(\hat{\bb}^\bo_S) = \lambda \hat{\bz}^\bo_S,
\end{equation}
where $\hat{\bz}^\bo_S \in \partial \|\hat{\bb}^\bo_S\|_1$. Combining equation (\ref{eq:oracle-gradient}) and (\ref{eq:oracle-penalty}), we obtain
\begin{equation} \label{eq:oracle-zero-subgradient}
     \nabla \mL_{\alpha, n}(\hat{\bb}^\bo_S)  - \nabla q_{\lambda}(\hat{\bb}^\bo_S) + \lambda \hat{\bz}^\bo_S = \b0.
\end{equation}
Hence $\hat{\bb}^\bo_S$ satisfies the zero-subgradient condition on the restricted program (\ref{eq:restircted}).
 By step (i) $\hat{\bb}_S$ is an interior point of the program (\ref{eq:restircted}), then it must also satisfy the zero-subgradient condition on the restricted program. Using the strict convexity from Lemma \ref{Lemma: restricted-strict-convex}, we obtain $\hat{\bb}_S = \hat{\bb}^\bo_S$.
 
\QEDB\\

The following lemma guarantees that the program in (\ref{eq:restircted}) is strictly convex:
\begin{lemma} \label{Lemma: restricted-strict-convex}
Suppose $\mL_{\alpha, n}$ satisfies the uniform RSC condition (\ref{as:RSC}) and $\rho_{\lambda}$ is $\mu$-amenable. Suppose in addition the sample size satisfied $n > \frac{2\tau}{\gamma - \mu}s \log p$, then the restricted program in (\ref{eq:restircted}) is strictly convex.
\end{lemma}
We omit the proof since it is similar to the proof of Lemma 2 in \cite{loh2017support}.

\QEDB\\

\textbf{Proof of step (ii) }: 
We rewrite the zero-subgradient condition (\ref{eq:zero-subgradient}) as 
\begin{equation*}
    \left (\nabla \mL_{\alpha, n}(\hat{\bb}) - \nabla \mL_{\alpha, n}(\bb^*) \right ) + \left (\nabla \mL_{\alpha, n}(\bb^*)- \nabla q_{\lambda}(\hat{\bb}) \right )+\lambda \hat{\bz} = \b0.
\end{equation*}
Let $\Hat{Q}$ be a $p \times p$ matrix $\hat{Q}= \int_0^1 \nabla^2 \mL_{\alpha, n} \blp \bb^* + t(\hat{\bb} - \bb^*) \brp dt$.  By the zero-subgradient condition and the fundamental theorem of calculus, we have
\begin{equation*}
    \hat{Q}(\hat{\bb} - \bb^*) + \left ( \nabla \mL_{\alpha, n}(\bb^*) - \nabla q_{\lambda}(\hat{\bb}) \right ) + \lambda \hat{\bz}=\b0.
\end{equation*}
And its block form is
\begin{equation} \label{eq: block-zero-subgradient}
    \left[ {\begin{array}{*{20}c}
   \hat{Q}_{SS} & \hat{Q}_{SS^c}  \\
   \hat{Q}_{S^c S} & \hat{Q}_{S^c S^c}  \\    
 \end{array} } \right]
  \left[ {\begin{array}{*{20}c}
   \hat{\bb}_S - \bb^*_S \\
   \b0 \\    
 \end{array} } \right]
 + \left[ {\begin{array}{*{20}c}
   \nabla \mL_{\alpha, n}(\bb^*)_S - \nabla q_{\lambda}(\hat{\bb}_S) \\
   \nabla \mL_{\alpha, n}(\bb^*)_{S^c} - \nabla q_{\lambda}(\hat{\bb}_{S^c}) \\    
 \end{array} } \right] + \lambda \left[ {\begin{array}{*{20}c}
   \hat{\bz}_{S} \\
   \hat{\bz}_{S^c} \\    
 \end{array} } \right] = \b0.
\end{equation}

The selection property implies $\nabla q_{\lambda}(\hat{\bb}_{S^c})=\b0$. Plugging this result into equation (\ref{eq: block-zero-subgradient}) and performing some algebra, we conclude that
\begin{equation}
    \hat{\bz}_{S^c} = \frac{1}{\lambda} \left \{ \hat{Q}_{S^cS}(\bb^*_S - \hat{\bb}_S) - (\nabla \mL_{\alpha, n}(\bb^*))_{S^c} \right \}.
\end{equation}
Therefore,
\begin{equation} \label{ineq:pd-feasible}
\begin{array}{ll}
\|\hat{\bz}_{S^c}\|_\infty  & \le \frac{1}{\lambda}\|\hat{Q}_{S^c S}(\hat{\bb}_S - \bb^*_S)\|_\infty + \frac{1}{\lambda} \|\nabla \mL_{\alpha, n}(\bb^*))_{S^c}\|_\infty\\
& \le \frac{1}{\lambda} \left \{ \max_{j \in S^c} \| e^T_j  \hat{Q}_{S^c S} (\hat{\bb}_S - \bb^*_S)\|_2 \right \} + \frac{1}{\lambda} \|\nabla \mL_{\alpha, n}(\bb^*))_{S^c}\|_\infty\\
& \le \frac{1}{\lambda} \left \{\max_{j \in S^c} \| e^T_j  \hat{Q}_{S^c S}\|_2  \right \} \| (\hat{\bb}_S - \bb^*_S)\|_2 + \frac{1}{\lambda} \|\nabla \mL_{\alpha, n}(\bb^*))_{S^c}\|_\infty.
\end{array}
\end{equation}
Observe that 
\begin{equation*}
\begin{array}{ll} 
 [(e_j^T\hat{Q}_{S^cS})_m]^2 & \le [\frac{1}{n}\sum_{i=1}^n  w(\bx_i)\bx_{ij} v(\bx_i) \bx_{im}\int_0^1 l''((y_i - \bx_i^T \bb^* - t(\bx_i\hat{\bb} - \bx_i\bb^*))v(\bx_i))\diff t  ]^2\\
 & \le  k_2^2 [\frac{1}{n} \sum_{i=1}^n w(\bx_i)\bx_{ij} \cdot v(\bx_i) \bx_{im} ]^2,\\
\end{array}
\end{equation*}
for all $j \in S^c$ and $m \in S$, where the second inequality follows from Assumption \ref{as:loss}(ii). By condition of Theorem \ref{Tm: oracle}, the variables $w(\bx_i)\bx_{ij}$ and $v(\bx_i) \bx_{im}$ are both sub-Gaussian. Using standard concentration results for i.i.d sums of products of sub-Gaussian variables, we have
\begin{equation*}
     P([(e_j^T\hat{Q}_{S^cS})_m]^2 \le c_1) \ge 1-c_2\exp(-c_3n).
\end{equation*}
It then follows from union inequality that 
\begin{equation} \label{ineq:maxQ}
    P( \max_{j \in S^c} \| e^T_j  \hat{Q}_{S^c S}\|_2 \le \sqrt{c_1s}) \ge 1 - c_2\exp(-c_3n + \log(s(p-s))) \ge 1-c_2\exp(-\frac{c_3}{2}n),
\end{equation}
where $n \ge \frac{2}{c_3} \log \left(s(p-s) \right)$.
By Lemma \ref{Lemma: restrict-bound}  we obtain
\begin{equation} \label{ineq:pd-feasible2}
    \|\hat{\bb}^\bo_S-\bb^*_S\|_2 \le C_6s^{-1/2}(\frac{\log p}{n})^{\zeta}.
\end{equation}
Furthermore, by (\ref{eq:gradient_bound}) we have
\begin{equation} \label{ineq:pd-feasible3}
    \|\nabla \mL_{\alpha, n}(\bb^*))_{S^c}\|_\infty \le \|\nabla \mL_{\alpha, n}(\bb^*))\|_\infty \le C_5(\frac{\log p}{n})^\zeta,
\end{equation}
Combining inequality (\ref{ineq:pd-feasible}), (\ref{ineq:maxQ}), (\ref{ineq:pd-feasible2}) and (\ref{ineq:pd-feasible3}), we have
\begin{equation*}
    \|\hat{\bz}_{S^c}\|_\infty \le \frac{1}{\lambda}C_7(\frac{\log p}{n})^\zeta.
\end{equation*}
with probability at least $1-C_8 \exp(-C_{41}s^{-\frac{1}{\zeta}}\log p)$. Note that $\alpha$ is required to satisfy both ranges in Theorem \ref{Tm: est} and (\ref{ineq: alpha_range2}). Combing these two ranges we have 
\begin{equation*}
   C_{22}s^{\frac{1}{k-1}}(\frac{n}{\log p})^{\frac{\zeta}{k-1}} \le \alpha \le C_{31} (\frac{n}{\log p})^{1-\zeta},
\end{equation*}
where $s = \bo \left ((\frac{n}{\log p})^{k(1-\zeta)-1} \right)$.
In particular, for $\lambda =C_{71}(\frac{n}{\log p})^{\zeta}$ with $C_{71}>C_7$, we conclude at last that the strict dual feasibility condition $ \|\hat{\bz}_{S^c}\|_\infty < 1$ holds, completing step (ii) of the PDW construction.

\textbf{Proof of step (iii) }: Since the proof for this step is almost identical to the proof in Step (iii) of Theorem 2 in \cite{loh2017statistical}, except for the slightly different notation, we refer the reader to  the arguments provided in that paper. 

\QEDB\\

To prove Theorem \ref{Tm: normality}, we need to generalized the asymptotic normality results for lower dimensional non-penalized M-estimator from \cite{he2000parameters} to the following Lemma:
\begin{lemma} \label{Lemma:Shao2000}
Suppose $\bz_1, \bz_2, \dots, \bz_n \in \RR^{p}$ are independent observations from probability distribution $F_{i,\bb}$, $i=1, 2, \dots, n$, with a common parameter $\bb \in \RR^s$. And $s$ may increase with the sample size $n$. Suppose $\mL_n(\bb) = \frac{1}{n}\sum_{i=1}^{n} \rho(\bz_i, \bb)$ is convex in $\bb$ in a neighborhood of $\bb^*$ and has a unique local minimizer $\hat{\bb}$. Define $\psi(\bz_i,\bb)= \frac{\partial}{\partial \bb} \rho(\bz_i, \bb) $ and $\eta_i(\tilde{\bb}, \bb) = \psi(\bz_i, \tilde{\bb}) - \psi(\bz_i, \bb) -E\psi(\bz_i, \tilde{\bb}) + E\psi(\bz_i, \bb) $ and $B_s = \{ \bnu \in \RR^s: \|\bnu\|_2=1\}$. Suppose $\bb^* \in \RR^s$ such that 
\begin{equation} \label{eq: true_beta}
    \|\sum_{i=1}^n \psi(\bz_i, \bb^*)\|_2 = \bo_p((ns)^{1/2}).
\end{equation}
Assume the following conditions are satisfied:
\begin{enumerate}
    \item[(i)] $\|\sum_{i=1}^n \psi(\bz_i, \hat{\bb})\|_2 = o_p(n^{1/2})$.
    \item[(ii)] There exist $C$ and $r \in (0,2]$ such that $\max_{i \le n} E_{\bb} \sup_{\tilde{\bb}:\|\tilde{\bb}-\bb\|_2 \le d} \|\eta_i(\tilde{\bb}, \bb)\|_2^2 \le n^C d^r$, for $0 <d \le 1$.
    \item[(iii)] There exists a sequence of $s$ by $s$ matrices $\bD_n$ with $\lim \inf_{n \rightarrow \infty} \lambda_{min}(\bD_n) > 0$ such that for any $K>0$ and uniformly in $\bnu \in B_s$,
    \begin{equation*}
        \sup_{\|\bb-\bb^*\|_2 \le K(s/n)^{1/2}} |\bnu^T\sum_{i=1}^{n} E_{\bb^*}(\psi(\bz_i, \bb) - \psi(\bz_i, \bb^*)) - n \bnu^T\bD_n(\bb - \bb^*)|=o((ns)^{1/2}).
    \end{equation*}
    \item[(iv)] $\sup_{\tilde{\bb}: \|\tilde{\bb} - \bb\|_2 \le K(s/n)^{1/2}} \sum_{i=1}^{n} E_{\bb} |\bnu^T \eta_i(\tilde{\bb}, \bb)|^2 =\bo(A(n,s))$ for any $\bb \in \RR^s$, $\bnu \in B_s$ and $K>0$.
    \item[(v)] $\sup_{\bnu \in S_s}\sup_{\tilde{\bb}: \|\tilde{\bb} - \bb\|_2 \le K(s/n)^{1/2}} \sum_{i=1}^{n} (\bnu^T \eta_i(\tilde{\bb}, \bb))^2 =\bo_p(A(n,s))$ for any $\bb \in \RR^s$ and $K>0$.  
\end{enumerate}
If $A(n,s)=o(n/\log n)$, we have 
\begin{equation*}
    \|\hat{\bb} - \bb^*\|_2^2=\bo_p(s/n).
\end{equation*}
Furthermore, if $A(n,s)=o(n/(s \log n))$, then for any unit vector $\bnu \in \RR^s$, we have
\begin{equation*}
    \hat{\bb} - \bb^* = - n^{-1} \sum_{i=1}^n \bD_n^{-1} \psi(\bz_i, \bb^*) +\br_n,
\end{equation*}
with $\|\br_n\|_2 = o_p(n^{-1/2})$.
\end{lemma}
\textit{Proof.} Our proof is similar to the proof of Theorem 1 and 2 in \cite{he2000parameters}. Note that in that paper,  $\bb^*$ is defined to be the solution of $\sum_{i=1}^n E_{\bb} \psi(x_i, \bb) = \b0$, in addition to the condition in equation (\ref{eq: true_beta}). However, a careful inspection of the proofs in that paper reveals that the results still holds if we only require $\bb^*$ to satisfied equation (\ref{eq: true_beta}).

\QEDB\\

{\bf Proof of Theorem \ref{Tm: normality}}\\
We then apply the result to the oracle estimator $\Hat{\bb}_S^{\bo}$ defined in equation (\ref{eq:oracle}) with $w(\bx) \equiv v(\bx)\equiv 1$. Although Lemma \ref{Lemma:Shao2000} requires $\mL_n$ to be convex, a throughout examination of the proofs in \cite{he2000parameters} shows that the results still hold if we restrict our attention to a subset of $\RR^p$ on which $\mL_n$ is convex and $\hat{\bb}$ is the unique minimizer. Since $\Hat{\bb}_S^{\bo}$ is $s$-dimensional vector without sparsity, we denote $\bx_i$, $\bb$ and $\bb^*$ all as $s$-dimensional vectors for the rest of our proof. We also denote $\bb_{\alpha}^*$ as $(\bb_{\alpha}^*)_S$. Let $\bz_i=(\bx_i, y_i)$ and we rewrite $\rho(\bz_i, \bb)$ as $l_\alpha(y_i - \bx_i^T\bb)$, with $\mL_{\alpha, n}$ taking the place of $\mL_{n}$. Then $\psi(\bz_i, \bb)= -l_\alpha'(y_i - \bx_i^T\bb)\bx_i$.

We start with verifying equation (\ref{eq: true_beta}), which can be rewrited as
\begin{equation} \label{eq: true_beta2}
     \|\sum_{i=1}^n l_\alpha'(\epsilon_i)\bx_i\|_2 = \bo_p((ns)^{1/2}).
\end{equation}
Observe that for any $\bnu \in B_s$,
\begin{equation} 
\begin{array}{ll} 
    P(|\sum_{i=1}^n \bnu^Tl_\alpha'(\epsilon_i)\bx_i - \sum_{i=1}^n E[\bnu^Tl_\alpha'(\epsilon_i)\bx_i]|>t) &\le nVar(\bnu^Tl_\alpha'(\epsilon_i)\bx_i)t^{-2} \\
    &\le nE|\bnu^Tl_\alpha'(\epsilon_i)\bx_i|^2t^{-2} \\
    &\le n E\|l_\alpha'(\epsilon_i)\bx_i\|_2^2t^{-2}\\
    & \le nsd_4t^{-2},
    \end{array}
\end{equation}
where the last inequality follows from inequality (\ref{ineq:psi_squred}). We then have
\begin{equation} \label{ineq:Opbound}
    P(|\sum_{i=1}^n \bnu^Tl_\alpha'(\epsilon_i)\bx_i|>t + \sum_{i=1}^n |E [\bnu^Tl_\alpha'(\epsilon_i)\bx_i]|)\le nsd_4t^{-2}.
\end{equation}
Observe that 
\begin{equation} \label{ineq: bound E-norm}
  \begin{array}{ll}  
 |E [\bnu^Tl_\alpha'(\epsilon_i)\bx_i]| &= |E[l_\alpha'(y_i - \bx_i^T\bb^*)\bnu^T\bx_i ]|\\
  &=|E[l_\alpha'(y_i - \bx_i^T\bb^*)\bnu^T\bx_i ]-E[l_{\alpha}'(y_i-\bx_i^T\bb_{\alpha}^*)\bnu^T\bx_i ]|\\
&  \le k_2E[|\bx_i^T(\bb_{\alpha}^*-\bb^*)||\bnu^T\bx_i|]\\
& \le  k_2\{E|\bx_i^T(\bb_{\alpha}^*-\bb^*)|^2\}^{\frac{1}{2}}\{E|\bnu^T\bx_i|^2\}^{\frac{1}{2}}\\
 & \le k_0^2k_2\|\bb_{\alpha}^*-\bb^*\|_2,
 \end{array}
\end{equation}
where the first and last inequalities follow from Assumption \ref{as:loss}(ii) and Assumption \ref{as:x}(iv) respectively. Together with the results in Theorem \ref{Tm:1} and condition $\alpha ^{1-k} = o(n^{-1/2})$, we obtain
\begin{equation} \label{eq:E-norm-bound}
    E[\bnu^Tl_\alpha'(\epsilon_i)\bx_i]= o(n^{-1/2}).
\end{equation}
Thus by inequality (\ref{ineq:Opbound}) and (\ref{eq:E-norm-bound}) we have $\sum_{i=1}^n \bnu^Tl_\alpha'(\epsilon_i)\bx_i = \bo_p((ns)^{1/2})$ for any $\bnu \in B_s$. It then implies equation (\ref{eq: true_beta2}).
Next we verify the conditions (i)-(v). Since the $\mL_{\alpha,n}$ is differentiable, the left hand side of condition (i) is 0 and thus it is satisfied. By definition of $\eta_i$, we have
\begin{equation*}
    \eta_i(\tilde{\bb}, \bb) = l_\alpha'(y_i-\bx_i^T\bb)\bx_i - l_\alpha'(y_i - \bx_i^T\tilde{\bb})\bx_i  -  E l_\alpha'( y_i - \bx_i^T\bb)\bx_i +
    E l_\alpha'(y_i - \bx_i^T\tilde{\bb})\bx_i
   . 
\end{equation*}
Observe that 
\begin{equation*}
\begin{array}{ll}
    \| \eta_i(\tilde{\bb}, \bb)\|_2 &\le \|l_\alpha'(y_i - \bx_i^T\tilde{\bb})\bx_i  - l_\alpha'(y_i-\bx_i^T\bb)\bx_i\|_2 + \| E l_\alpha'(y_i - \bx_i^T\tilde{\bb})\bx_i -
    E l_\alpha'(y_i-\bx_i^T\bb)\bx_i\|_2\\
    &  \le k_2|\bx_i^T(\tilde{\bb} - \bb)|\cdot\|\bx_i\|_2 + k_2\|E \bx_i^T(\tilde{\bb} - \bb)\bx_i\|_2 \\
    & \le  k_2\|\tilde{\bb} - \bb\|_2 \|\bx_i\|_2^2 +  k_2E\|\bx_i^T(\tilde{\bb} - \bb)\bx_i\|_2\\
    & \le  k_2\|\tilde{\bb} - \bb\|_2 \|\bx_i\|_2^2 +  k_2\|\tilde{\bb} - \bb\|_2 E \|\bx_i\|_2^2,
    \end{array}
\end{equation*}
where the second and third inequality follow from Assumption \ref{as:loss}(ii) and Jensen's inequality, respectively. 
We then obtain
\begin{equation*}
   \max_{i \le n} E_{\bb} \sup_{\tilde{\bb}:\|\tilde{\bb}-\bb\|_2 \le d} \|\eta_i(\tilde{\bb}, \bb)\|_2^2 \le \max_{i \le n} 4 k_2^2d^2 E \|\bx_i\|_2^4.
\end{equation*}
Since Assumption \ref{as:x}(iv) implies $E \|\bx_i\|_2^4 =\bo (s^2)$ for $i=1,\cdots,n$, condition(ii) holds with $r=2$ and if $s=\bo(n^{r_1})$ for $r_1>0$.

Similarly, for any $\bnu \in B_s$, we have
\begin{equation*}
\begin{array}{ll}
     \bnu^T \eta_i(\tilde{\bb}, \bb) &\le |l_\alpha'( y_i - \bx_i^T\tilde{\bb}) -l_\alpha'(y_i-\bx_i^T\bb)||\bnu^T\bx_i| +  E [|l_\alpha'( y_i - \bx_i^T\tilde{\bb}) -l_\alpha'(y_i-\bx_i^T\bb)||\bnu^T\bx_i|]\\
    &  \le  k_2|\bx_i^T(\tilde{\bb} - \bb)||\bnu^T \bx_i| +  k_2E[|\bx_i^T(\tilde{\bb} - \bb)||\bnu^T \bx_i|]\\
    & \le  k_2\|\tilde{\bb} - \bb\|_2 |\tilde{\bnu}^T\bx_i|| \bnu^T \bx_i | +  k_2\|\tilde{\bb} - \bb\|_2 \{E|\tilde{\bnu}^T\bx_i|^2\}^{1/2} E\{| \bnu^T \bx_i |^2\}^{1/2}\\
    & \le  k_2\|\tilde{\bb} - \bb\|_2(|\tilde{\bnu}^T\bx_i|| \bnu^T \bx_i | + k_0^2),
    \end{array}
\end{equation*}
where $\tilde{\bnu}=(\tilde{\bb} - \bb)/\|\tilde{\bb} - \bb\|_2$. The second and last inequalities follow from Assumption \ref{as:loss}(ii) and Assumption \ref{as:x}(iv) respectively. It then gives
\begin{equation*}
    | \bnu^T \eta_i(\tilde{\bb}, \bb)|^2 \le k_2^2\|\tilde{\bb} - \bb\|_2^2(|\tilde{\bnu}^T\bx_i|^2| \bnu^T \bx_i |^2 + 2k_0^2|\tilde{\bnu}^T\bx_i|| \bnu^T \bx_i |+k_0^4).
\end{equation*}
Together with Assumption \ref{as:x}(iv), we obtain
\begin{equation}
    E| \bnu^T \eta_i(\tilde{\bb}, \bb)|^2 = \bo(\|\tilde{\bb} - \bb\|_2^2)
\end{equation}
and 
\begin{equation}
    | \bnu^T \eta_i(\tilde{\bb}, \bb)|^2 = \bo_p(\|\tilde{\bb} - \bb\|_2^2).
\end{equation}
Hence condition (iv) and (v) hold with $A(n,s)=s$.

Finally we show condition (iii). Let $\bD_{\alpha, n} = E[\nabla^2 \mL_{\alpha, n}(\bb^*)]$ and thus it is an $s$ by $s$ matrix. Observe that
\begin{equation*}
\begin{array}{ll}
      E[l_{\alpha}''(\epsilon_i)|\bx_i] & = E[l''_\alpha (\epsilon_i)\mathbbm{1}(|\epsilon_i| \le \alpha)|\bx_i ] + E[l''_\alpha (\epsilon_i)\mathbbm{1}(|\epsilon_i| > \alpha)|\bx_i] \\
      &  \ge E[(1-d_1|\epsilon_i|^k\alpha^{-k}) \mathbbm{1}(|\epsilon_i| \le \alpha)|\bx_i ] + E[l''_\alpha (\epsilon_i)\mathbbm{1}(|\epsilon_i| > \alpha)|\bx_i]\\
     & \ge  P(|\epsilon_i| \le \alpha |\bx_i) - d_1\alpha^{-k}E[|\epsilon_i|^k|\bx_i]- k_2\alpha^{-k}E[|\epsilon_i|^k|\bx_i]
      \\
     & \ge 1-\alpha^{-k}E[|\epsilon_i|^k|\bx_i] - d_1\alpha^{-k}E[|\epsilon_i|^k|\bx_i]- k_2\alpha^{-k}E[|\epsilon_i|^k|\bx_i]\\
     &= 1- (d_1+k_2+1)\alpha^{-k}E[|\epsilon_i|^k|\bx_i],
\end{array}
\end{equation*}
where the first and second inequalities follow from Assumption \ref{as:loss}(iii) and (ii), respectively. Thus for any $\bnu \in B_s$, we have
\begin{equation*}
    \begin{array}{ll}
    \bnu^T \bD_{\alpha, n} \bnu & = E[l_{\alpha}''(\epsilon_i) \bnu^T\bx_i \bx_i^T \bnu]\\
     & \ge E[(1- (d_1+k_2+1)\alpha^{-k} E[|\epsilon_i|^k|\bx_i])\bnu^T\bx_i \bx_i^T \bnu]
     \\
    & = \bnu^T E[\bx_i \bx_i^T] \bnu - (d_1+k_2+1) \alpha^{-k} E[E(|\epsilon_i|^k|\bx_i) (\bnu\bx_i)^2 ] \\
     &  \ge k_l - (d_1+k_2+1) \alpha^{-k} \{ E[E(|\epsilon_i|^k|\bx_i)]^2 \}^{1/2} \{E[(\bnu\bx_i)^4]\}^{1/2} \\
     & \ge k_l  - C_9\alpha^{-k},
    \end{array}
\end{equation*}
where second inequality follows from Assumption \ref{as:x}(i) and  $C_9$ is a constant that only depends on $k_0$, $k_2$, $d_1$, $M_k$. Hence if $\alpha > (2C_9/k_l)^{1/k}$, we have $\lambda_{min}( \bD_{\alpha, n}) >k_l/2$. It then implies $\lim \inf_{n \rightarrow \infty} \lambda_{min}(\bD_{\alpha, n}) > 0$. Observe that 
    \begin{equation*}
    \begin{array}{ll}
        &|\bnu^T\sum_{i=1}^{n} E_{\bb^*}(\psi(\bx_i, \bb) - \psi(\bx_i, \bb^*)) - n \bnu^T\bD_{\alpha, n}(\bb - \bb^*)|\\
        = & |\bnu^T\sum_{i=1}^{n} E_{\bb^*}\{(l'_\alpha( y_i - \bx_i^T\bb^*)\bx_i - l'_\alpha(y_i - \bx_i^T\bb)\bx_i - l''_\alpha(y_i - \bx_i^T\bb^*)\bx_i\bx_i^T(\bb - \bb^*)\}|\\
         = & |\bnu^T\sum_{i=1}^{n} E_{\bb^*}\{(l''_\alpha(y_i - \bx_i^T\tilde{\bb})\bx_i^T(\bb - \bb^*)\bx_i - l''_\alpha( y_i - \bx_i^T\bb^*)\bx_i\bx_i^T(\bb - \bb^*)\}|\\
       \le  & |\bnu^T\sum_{i=1}^{n} E_{\bb^*}\{(k_3|\bx_i^T(\tilde{\bb} -\bb^*)| |\bx_i^T(\bb - \bb^*)\bx_i| \}| \\
       \le & k_3\|\bb - \bb^*\|_2^2\sum_{i=1}^{n} E_{\bb^*}\{|\bx_i^T\tilde{\bnu}|^2 |\bx_i^T\bnu| \},
    \end{array}
    \end{equation*}
where $\tilde{\bb}$ is a vector lying between $\bb$ and $\bb^*$ and $\tilde{\bnu}=(\tilde{\bb} - \bb)/\|\tilde{\bb} - \bb\|_2$. Note that the second equation follows from mean value theorem and the first inequality follows from the condition that $l_{\alpha}''$ is Lipschitz. By Assumption \ref{as:x} (iv) we have $\sum_{i=1}^{n} E_{\bb^*}\{|\bx_i^T\tilde{\bnu}|^2 \bx_i^T\bnu| \} = \bo(n)$. Hence condition (iii) holds if $s/n \rightarrow 0$.

We conclude that the desired results hold for the oracle estimator $\hat{\bb}^\bo_S$. In particular, we have 
\begin{equation} \label{eq:est-error-sum}
\begin{split}
     \hat{\bb}^\bo_S - \bb^* = & n^{-1} \sum_{i=1}^n \bD_{\alpha, n}^{-1} l_\alpha'(\epsilon_i)\bx_i + \br_n\\
     = &  n^{-1} \sum_{i=1}^n \{\bD_{\alpha, n}^{-1} l_\alpha'(\epsilon_i)\bx_i - E[\bD_{\alpha, n}^{-1} l_\alpha'(\epsilon_i)\bx_i] \} + E[\bD_{\alpha, n}^{-1} l_\alpha'(\epsilon_i)\bx_i]+ \br_n,
\end{split}
\end{equation}
with $\|\br_n\|_2 = o_p(n^{-1/2})$. Observe that
\begin{equation} \label{eq:mu-bound}
\begin{array}{ll}
     \| E[\bD_{\alpha, n}^{-1}l_\alpha'(\epsilon_i)\bx_i] \|_2 &= \| \bD_{\alpha, n}^{-1}E[l_\alpha'(\epsilon_i)\bx_i] \|_2\\
     &=  \|\bD_{\alpha, n}^{-1} \tilde{\bnu}\|_2\|E[l_\alpha'(\epsilon_i)\bx_i]\|_2\\
     & \le [\lambda_{\min}(\bD_{\alpha, n})]^{-1} \|E[l_\alpha'(\epsilon_i)\bx_i]\|_2 \\
     &= o(n^{-1/2}),
\end{array}
\end{equation}
where the last equality follows from equation (\ref{eq:E-norm-bound}). By equations (\ref{eq:est-error-sum}) and (\ref{eq:mu-bound}), we obtain
\begin{equation} \label{eq: asymptotic}
    \frac{\sqrt{n}}{\sigma_\bnu} \cdot \bnu^T(\hat{\bb}^\bo_S - \bb^*) \xrightarrow{d} N(0,1),
\end{equation}
where $\sigma_\bnu^2 = \bnu^T \bD_{\alpha, n}^{-1} Var(l_\alpha'(\epsilon_i)\bx_i)  \bD_{\alpha, n}^{-1} \bnu$. By Theorem \ref{Tm: oracle}, the asymptotic result in (\ref{eq: asymptotic}) is also applicable for the stationary point  $\tilde{\bb}$.

\QEDB\\
To prove Theorem \ref{Tm: RSC}, we need the following result:
\begin{lemma}\label{Lemma: RSC_loh}
Suppose covariate $\bx$ satisfies Assumption 3(iv) and $\ell''_\alpha(u)$ satisfies Assumption 2(ii). For any fixed $\alpha>0$, suppose  the bound $ C_{14}k_0^2 \left(\sqrt{\varepsilon_{T}}+\exp\left(-\frac{C_{15}T^2}{k_0^2r^2}\right)\right) < \frac{\gamma_{\alpha, T}}{\gamma_{\alpha, T} + k_2} \cdot \frac{k_l}{2}$ holds, where $ \gamma_{\alpha, T} = \min_{|u| \le T}\ell''_{\alpha}(u)>0$. Suppose in addition that the sample size satisfies $n \ge C_{10}s \log p$. With probability at least $1-C_{11}\exp(-C_{12}\log p)$, the loss function $\mL_{\alpha,n}$ satisfies 
\begin{equation} \label{eq:RSC}
    	\langle \nabla \mL_{\alpha,n}(\bb_1)- \nabla \mL_{\alpha,n}(\bb_2), \bb_1-\bb_2 \rangle \ge \gamma_{\alpha} \|\bb_1-\bb_2\|_2^2 - \tau_{\alpha} \frac{\log p}{n}\|\bb_1 - \bb_2\|_1^2,
\end{equation}
where $\bb_j \in \RR^p$ such that $\|\bb_j-\bb^*\|_2\le r$ for $j=1, 2$ with 
\begin{equation} \label{eq:RSC_coeff}
    \gamma_{\alpha}=\frac{\gamma_{\alpha, T}k_l}{16} \quad \text{and} \quad \tau_{\alpha}=\frac{C_{13}(\gamma_{\alpha, T} + k_2)^2k_0^2T^2}{r^2}.
\end{equation}
Here the constants $C_{10}, C_{11}, C_{12}, C_{13}, C_{14}, C_{15}$ do not depend on $\alpha$.
\end{lemma}
Proof. The proof is similar to the proof of Proposition 2 in \cite{loh2017statistical}. Note that in that paper, it assumes $\bx_i \indep \epsilon_i$. However, a careful inspection of the proofs reveals that the result stills holds if we allow $\epsilon_i$ to depend on $\bx_i$. We refer the reader to the arguments provided in that paper.

\QEDB\\

{\bf Proof of Theorem \ref{Tm: RSC}}\\
Recall $\gamma_{\alpha, T} = \min_{|u| \le T} \ell''_{\alpha}(u)$. By Assumption 2(iii) , $\alpha \ge \alpha_0$ and $\alpha_0=\max\{ (2d_1)^{\frac{1}{k}},1\}\cdot T_0$ we have

\begin{equation} \label{eq:lb_gamma}
    \gamma_{\alpha, T_0} = \min_{|u| \le T_0} \ell_{\alpha}(u)
    \ge  \min_{|u| \le T_0} (1-d_1|u|^{k}\alpha^{-k}) \ge 1-d_1|T_0|^{k}\alpha_0^{-k} \ge \frac{1}{2}.
\end{equation}
And 
\begin{equation} \label{eq:ub_gamma}
    \gamma_{\alpha, T_0} = \min_{|u| \le T_0} \ell_{\alpha}(u)
    \le  \min_{|u| \le T_0} (1+d_1|u|^{k}\alpha^{-k}) \le 1+d_1|T_0|^{k}\alpha_0^{-k} \le \frac{3}{2}.
\end{equation}
By equation (\ref{eq:lb_gamma}), we obtain
\begin{equation*}
\frac{\gamma_{\alpha, T_0}}{\gamma_{\alpha, T_0} + k_2} \cdot \frac{k_l}{2}  \ge \frac{\frac{1}{2}}{\frac{1}{2} + k_2} \cdot \frac{k_l}{2} \ge \frac{k_l}{2+4k_2}.
\end{equation*}
Together with condition $C_{14}k_0^2 \left(\sqrt{\varepsilon_{T_0}}+\exp \left(-\frac{C_{15}T_0^2}{k_0^2r^2} \right) \right) < \frac{k_l}{2+4k_2}$, we have
\begin{equation} \label{eq:RSC_condition}
    c_{14}k_0^2 \left(\sqrt{\varepsilon_{T_0}}+\exp \left(-\frac{c_{15}T_0^2}{k_0^2r^2} \right) \right) < \frac{\gamma_{\alpha, T_0}}{\gamma_{\alpha, T_0} + k_2} \cdot \frac{k_l}{2}.
\end{equation}
By equation (\ref{eq:lb_gamma}), (\ref{eq:ub_gamma}), (\ref{eq:RSC_condition}) and Lemma \ref{Lemma: RSC_loh} we complete the proof.

\QEDB\\
\bibliography{pram}

\end{document}